\documentclass[a4paper, 11pt]{article}
\usepackage{setspace}
\pdfoutput=1

\usepackage{makeidx}
\usepackage{amssymb}
\usepackage{amsmath}
\usepackage{amsfonts}
\usepackage{lineno}
\usepackage{epsfig}
\usepackage{graphics}
\usepackage{graphicx}
\usepackage{ntheorem}
\usepackage{color}
\usepackage[top=1in,bottom=1in,left=1.2in,right=1.2in]{geometry}
\usepackage{mathrsfs}
\usepackage{enumerate}
\usepackage[linesnumbered,ruled,vlined]{algorithm2e}
\usepackage{float}
\usepackage{subfigure}
\usepackage{caption}
\usepackage{stmaryrd}
\usepackage{epsfig}
\usepackage{xy}\xyoption{all}

\usepackage[bookmarks=false]{hyperref}

\newtheorem{theorem}{Theorem}
\newtheorem{assumption}{Assumption}

\newtheorem{definition}{Definition}

\newtheorem{lemma}{Lemma}

\newtheorem{proposition}{Proposition}
\newtheorem{remark}{Remark}

\newcommand{\te}{\hfill $\Box$}

\newcommand{\be}{\begin{equation}}
\newcommand{\ee}{\end{equation}}
\newcommand{\ba}{\begin{eqnarray}}
\newcommand{\ea}{\end{eqnarray}}
\newcommand{\bas}{\begin{eqnarray*}}
\newcommand{\eas}{\end{eqnarray*}}
\def\text{\hbox}

\DeclareMathAlphabet{\mathbfit}{OML}{cmr}{bx}{it}

\begin{document}

\title{Nonstationary Discounted Stochastic Games under Prospect Theory with Applications to the Smart Grid}
%\title{Discounted Stochastic Games for the Smart Grid with Prospect Prosumers}

%\author{Xi-Ren~Cao,~\IEEEmembership{Fellow,~IEEE}
\author{Yiting~Wu, Junyu~Zhang
\thanks{Y.T.~Wu is with the School of Mathematics, Sun Yat-Sen University, Guangzhou 510275, China (email: wuyt35@mail2.sysu.edu.cn). J.Y.~Zhang is with the School of Mathematics, Sun Yat-Sen University, Guangzhou 510275, China (corresponding author, email: mcszhjy@mail.sysu.edu.cn).}
}

%
%
%\markboth{IEEE TRANSACTIONS ON AUTOMATIC CONTROL,~Vol.~xx, No.~xx, July~2015}
%{Shell \MakeLowercase{\textit{et al.}}: Bare Demo of IEEEtran.cls for Journals}
\date{}
\maketitle

\begin{abstract}
	This paper considers the discounted criterion of nonzero-sum decentralized stochastic games with prospect players. The state and action spaces are finite. The state transition probability is nonstationary. Each player independently controls their own Markov chain.
	The subjective behavior of players is described by the prospect theory (PT). Compared to the average criterion of stochastic games under PT studied firstly in 2018, we are concerned with the time value of utility, i.e., the utility should be discounted in the future. Since PT distorts the probability, the optimality equation that plays a significant role in proving the existence of equilibrium does not exist. On the other hand, the games change into Markov decision processes (MDPs) with nonstationary payoff function when fixing others' stationary Markov strategies, then the occupation measure and the linear programming of stationary MDPs are no longer suitable. Therefore, we explore a new technique by constructing the marginal distribution on the state-action pairs at any time, and establish the existence of Nash equilibrium. When the initial state is given, there exists a Markov Nash equilibrium. Furthermore, this novel technique can be extended to the finite horizon criterion. Then, we present an algorithm to find a Markov $\varepsilon$-equilibrium. Finally, the model is applied to a noncooperative stochastic game among prosumers who can produce and consume energy in the smart grid, and we give some simulation results.
%\vspace{6cm}

\begin{center}
\end{center}

%Nonzero-sum stochastic games, expected discounted payoff criterion, Borel state space, unbounded payoffs, almost Markov $\varepsilon$-equilibrium

\end{abstract}

%
%\begin{IEEEkeywords}
%Confluencity, performance potential, bias potential, HJB equation,
%bias optimality, weak ergodicity, weak recurrence, direct-comparison
%based optimization
%\end{IEEEkeywords}
%
%\IEEEpeerreviewmaketitle
\begin{spacing}{1.0}
\end{spacing}
%\newpage

\section{Introduction}

According to whether players independently control the Markov chain, stochastic games can be classified as below: 1) centralized stochastic games (all players jointly control a single Markov chain); 2) decentralized stochastic games (each player independently controls their own Markov chain). There are numerous studies on centralized stochastic games, such as \cite{dufour,Huang2020,jaskiewicz2016,nowak,sennott,wei2019,wei2022}. 
However, the decentralized stochastic games mainly focus on the average constrained games. The primary techniques are the occupation measure and linear programming. Altman et al \cite{altman} introduces the discrete-time finite state and action spaces model, which they call cost-coupled constrained stochastic games with independent state processes. They establish the stationary constrained Nash equilibrium. Based on \cite{altman}, Singh and Hemachandra \cite{singh} characterize the Nash equilibrium through the global minimizers of a non-convex mathematical program. Zhang \cite{zhang2019} extends \cite{altman} to a countable state space and unbounded cost functions. Moreover, Etesami \cite{etesami_learning} proposes a learning algorithm of stationary $\varepsilon$-equilibrium% in decentralized stochastic games with long-term expected average criterion
, and applies to energy management in the smart grid.
For continuous-time model, Zhang and Zou \cite{zhang2022} study unbounded transition rates as well as cost rates and prove the existence of constrained Nash equilibrium. In this paper, we pay attention to decentralized stochastic games where players couple with each other through a payoff function, but we are not limited to the average constrained case.

Players, as decision makers in stochastic games, their behavior cannot be ignored. Saad et al \cite{Saadproc} discusses the subjective view of users about their opponents, and illustrates its impact on the outcome of energy management. The subjective behavior of players can be modelled by prospect theory (PT), a Nobel Prize winning theory \cite{kahneman}. PT emphasizes the reaction of players to risk and uncertainty, which changes the effect value by the valuation function and distorts the probability by the weighting function.
Etesami et al \cite{etesami} first studies a nonzero-sum decentralized stochastic game of players with subjective behavior, which considers the expected average criterion under PT. The expected average criterion represents that players focus on the long term and ignore the short term, but players are interested in the time value of the expected utility sometimes. Thus, the discounted criterion is extremely important, which players believe that the utility in future needs to be discounted (cannot be achieved immediately). We mainly consider the expected discounted criterion that is not studied for stochastic games under PT to the best of our knowledge.

For the general nonzero-sum stochastic games (without PT), the optimality equation is the crucial technique. There are different optimality equation forms of different criteria, such as expected discounted criterion \cite{sennott}, expected average criterion \cite{nowak}, finite horizon criterion \cite{wei2019} and probability criterion \cite{Huang2020}. All of them transform the stochastic game into a static game using the optimality equation, and come back to the original game by dynamic programming and so on. For the nonzero-sum stochastic games under PT, the state transition probability of the equivalent expected stochastic games is not Markov (depends on all history) due to the distortion of probability, thus the optimality equation does not exist. Therefore, the nonzero-sum stochastic games under PT need new tools.

In \cite{etesami}, Etesami et al proves the existence of stationary Nash equilibrium under certain conditions (refer to Remark \ref{remark1}). The model of stochastic game is stationary, including the state transition probability and the payoff function. It becomes a Markov decision process (MDP) of a player with stationary state transition probability and payoff function when fixing the stationary Markov strategies of other players. Then, Etesami et al transforms the stochastic game into a static game with the occupation measure and the linear programming of stationary MDP, thereby solving the stochastic game. However, fixing others' stationary Markov strategies, the discounted (or finite horizon) model (even if the state transition probability and the payoff function are stationary) changes into a
MDP with nonstationary payoff function. Though nonstationary MDP can be approximated by linear programming \cite{Ghate}, the method of \cite{etesami} or the primary techniques of decentralized stochastic games are no longer suitable. In this paper, we explore new techniques by constructing the marginal distribution on the state-action pairs at any time %(we call it the occupation measure with time) 
and the set that it takes value (i.e., (\ref{M})). Moreover, we extend the stochastic game of players under PT to the nonstationary case, and the payoff function also depends on state. Especially, the nonstationary decentralized stochastic game under the expected discounted criterion without PT is a special case of our model (see Remark \ref{thm_rem} (c)). It is important to note that most literature focuses on the stationary model, e.g., \cite{jaskiewicz2016,sennott,wei2022}. Nevertheless, we demonstrate the existence of Nash equilibrium for the nonstationary decentralized stochastic games without PT simultaneously.

The algorithm research of nonzero-sum stochastic games is a challenge, unlike MDP and zero-sum stochastic games, which have good properties, such as the contraction of operators and the uniqueness of the value function. Rao et al \cite{rao} develops two convergent algorithms for two-person zero-sum discounted stochastic games with finite state and action spaces, and considers a possible extension to nonzero-sum stochastic games. However, the monotone value sequence cannot be established for the nonzero-sum case. Since it is difficult to discover a general algorithm for nonzero-sum stochastic games, some works focus on the specific model. For instance, Raghavan and Syed \cite{raghavan} assume that the transition probability depends on only one of the players, and formulate a linear complementary problem to solve the Nash equilibrium. In addition, several existing algorithm researches on nonzero-sum stochastic games lack an analysis of convergence. Innovatively, \cite{etesami} presents an algorithm for obtaining equilibrium and proves its convergence. Now, we study the equilibrium algorithm for the nonzero-sum discounted stochastic games with PT under a sampling without replacement, inspired by \cite{etesami} which is a sampling with replacement method.

As an interdisciplinary field, stochastic game has wide applications in investment and reinsurance \cite{zhang_insurance}, energy management\cite{etesami}, resource management\cite{kaitala}, queuing theory\cite{nowak_book}, etc. In this work, we illustrate the applications of nonstationary discounted stochastic games under PT in the smart grid energy management. 

Differs from traditional centralized power generation, distributed energy networks not only focus on power demands but also on ecological benefits, e.g., low carbon emissions, low cost \cite{li2018}. Smart grid, as an important product of the development of the electric power grid in recent years, is equipped with renewable resource generation, storage devices, smart meters and so on \cite{atzeni}. At the moment, the energy management of distributed energy networks is an urgently significant research content.

Based on the equivalence relationship between prosumers who can produce energy by renewable resources and consume energy in the smart grid \cite{liu2020}, it forms a noncooperative game among prosumers. 
Mohsenian-Rad et al \cite{Mohsenian-Rad} first constructs demand response as a noncooperative game, and derives the optimal load schedule by Nash equilibrium. Each user aims to minimize the energy cost which depends on other users' load. Furthermore, Soliman and Leon-Garcia \cite{Soliman} consider the game among users with energy storage capacity. The optimization variables are energy storage, power consumption from the grid and the storage device.
The above literature study static game, and do not consider renewable resources. Actually, the renewable resource generation brings uncertainty due to the weather condition, which makes the stochastic game model more suitable. 
Etesami et al \cite{etesami} first investigates a stochastic game of prosumers who own renewable energy sources, managing the uncertainty by storage devices. 
Chen et al studies Peer-to-Peer energy sharing markets by stochastic games \cite{chen_data} and stochastic leader-follower games \cite{chen_model}. Though it lacks storage devices, the prosumers can manage energy by selling the redundant energy. Moreover, Chen et al \cite{chen_voltage} also discusses distributed voltage regulation through stochastic games. Etesami et al \cite{etesami_repeat} formulates the routing of electric vehicles as a repeated game, a special type of stochastic game.
However, the work about stochastic games in the smart grid are not very much. This paper is more attentive towards it. 

The main contribution of this paper is a novel nonstationary decentralized stochastic game model considering the discounted criterion under PT, which are as follows.
\begin{enumerate}[1)]
	\item For players with subjective behaviors, this paper first studies the discounted criterion of nonzero-sum discrete-time decentralized stochastic games with finite state and action space under PT.
	\item The considered model can be extended to the nonstationary case, i.e., the state transition probability and payoff function are nonstationary, which makes it more general.
	\item It proposes a new technique for the existence of Nash equilibrium by using the marginal distribution on the state-action pairs at any time. And this technique is appropriate in the finite horizon criterion with PT, see Remark \ref{thm_rem} (b).
	\item It provides a convergent algorithm of $\varepsilon$-equilibrium. % which \Red{compare to the average criterion of stochastic games under PT \cite{etesami}.} 
	As an application, we consider a novel and practical performance criterion for the stochastic games of prospect prosumers in the smart grid energy management, and analyzes the simulation results of $\varepsilon$-equilibrium algorithm.
\end{enumerate}

This paper is organized as follows. 
In Section \ref{modelsection}, we introduce the nonstationary discounted decentralized stochastic game model under PT. In Section \ref{proofsection}, we study the existence of Nash equilibrium for the stochastic games with prospect players. An algorithm of $\varepsilon$-equilibrium is provided in
Section \ref{algsection}. In Section \ref{applisection}, we give an application in the smart grid energy management. Finally, we conclude this paper in Section \ref{conclusection}.

\section{The Model of Stochastic Games}\label{modelsection}
\noindent \textbf{Notation}: For any topological space $X$, we endow it with Borel $\sigma$-algebra $\mathscr{B}(X)$, i.e., $\sigma$-algebra generated by all open subsets of $X$. $\mathscr{P}(X)$ denotes the set of all Borel probability measures on $X$ endowed with weak topology. $|X|$ denotes the number of elements in $X$.

In this section, we first introduce the nonzero-sum discounted decentralized stochastic games under PT. Denote the players set by $N:=\{1,2,\cdots,n\}$. We consider the discrete-time case, i.e., time set $T:=\{1,2,\cdots\}$. The model is as below.

\be \label{model}
%\mathbb{G}:=
\left\{\left(S_i,A_i,r_i(\boldsymbol{s},\boldsymbol{a}),q_t^i(\cdot|s^i,a^i) \right)_{i\in N, t\in T},\beta \right\},
\ee
consists of the following elements:
\begin{enumerate}
	\item $S_i$ is a nonempty finite state space of player $i$. Let $\boldsymbol{S}:=\prod_{i\in N}S_i$ be the state space of all players (or the system).
	\item $A_i$ is a nonempty finite action space of player $i$. Denote $\boldsymbol{A}:=\prod_{i\in N}A_i$, i.e., action space of all players.% $A_i(x)\subset A_i$ is the available action set for player $i\in I$ in state $x\in S$. 
	\item $r_i(\boldsymbol{s},\boldsymbol{a})$ is (real-valued) payoff function of player $i$, which depends on the current system state $\boldsymbol{s}=(s^1,s^2,\cdots,s^n)\in \boldsymbol{S}$ and the action $\boldsymbol{a}=(a^1,a^2,\cdots,a^n)\in \boldsymbol{A}$ taken by all players.\footnote{It is the key to forming a game, where players couple with each other through payoffs.}
	\item For player $i$, $q_t^i(\cdot|s^i,a^i)$ is state transition probability that describes the distribution of the next state based on state $s^i\in S_i$ and action $a^i\in A_i$ at time $t\in T$. Assume that the state transition probability of system $q{_t}(\cdot|\boldsymbol{s},\boldsymbol{a})$ satisfies
	\begin{align}\label{probabilityindependence}
		q_t(\boldsymbol{y}|\boldsymbol{s},\boldsymbol{a})=\prod_{i\in N}q_t^i(y^i|s^i,a^i),~~\forall t\in T,
	\end{align}
	for all $\boldsymbol{s}=(s^i:i\in N),\boldsymbol{y}=(y^i:i\in N)\in \boldsymbol{S},\boldsymbol{a}=(a^i:i\in N)\in \boldsymbol{A}$.\footnote{It is a symbol of decentralized stochastic games, i.e., each player independently controls their own Markov chain.}
	\item $\beta \in (0,1)$ is a discount factor.
\end{enumerate}

In stochastic game (\ref{model}), payoff function and action set can be related to time, as noted in Remark \ref{thm_rem} (c). We call this model a nonstationary stochastic game. For the sake of simplicity, we take the nonstationary state transition probability as an example to illustrate.

Let $H_1^i:=S_i, H_t^i:=(S_i\times A_i)^{t-1}\times S_i$ $(t\geq 2,t\in T)$ be the space of all admissible (state and action) histories of the game up to time $t$ for player $i$, and $H_{\infty}:=(\boldsymbol{S}\times \boldsymbol{A})^{\infty}$, $H_{\infty}^i:=(S_i\times A_i)^{\infty}$, all of which are endowed with product $\sigma$-algebra.

To define the optimality criterion, we introduce the definition of strategy for each player as below.

\begin{definition}\label{def1}
	\begin{enumerate}[(a)]
		\item A (randomized) history-dependent strategy for player $i\in N$ is a sequence $\pi^i=\{\pi^i_1,\pi^i_2,\cdots\}$ of stochastic kernels $\pi^i_t~(t\in T)$ on $A_i$ given $H_t^i$, i.e., $\pi^i_t(\cdot|{h}_t^i)\in \mathscr{P}(A_i)$ for any history ${h}_t^i=(s_1^i,a_1^i,s_2^i,a_2^i,\cdots,s_t^i)\in H_t^i$.%, and $\pi^i_t(B|\cdot)\in \mathscr{B}(H_t^i)$ for any $B\in \mathscr{B}(A_i)$. 
		%$H_t$, i.e., $\pi^i_t(\cdot|\{h_t})\in \mathscr{P}(A_i)$ for any history $\{h_t}=(x_1,\{a}_1,x_2,\{a}_2,\cdots,x_t)\in H_t$ with $\{a}_s=(a^i_s:i\in N)\in \{A}(s=1,2,\cdots,t-1)$, and $\pi^i_t(B|\cdot)\in \mathscr{B}(H_t)$ for any $B\in \mathscr{B}(A_i)$. 
		%	\item If $\pi_i^s\equiv \pi_i^t(\forall s,t\geq 1,s\neq t)$, then $\pi_i$ is called a stationary history-dependent policy for player $i$ and we denote $\pi_i=\{\pi_i,\pi_i,\cdots\}$.
		\item If there exists a sequence $\phi^i=\{\phi^i_1,\phi^i_2,\cdots\}$ of stochastic kernels $\phi^i_t~(t\in T)$ on $A_i$ given $S_i$ %satisfying $\phi^i_t(A_i|x)=1(\forall x\in S)$ 
		such that $\pi^i_t(\cdot|{h}_t^i)\equiv\phi_t^i(\cdot|s_t^i)~(\forall t\in T,{h}_t^i=(s_1^i,a_1^i,s_2^i,a_2^i,\cdots,s_t^i)\in H_t^i)$, then $\pi^i$ is called a Markov strategy for player $i$.
		\item If there exists a sequence $d^i=\{d^i_1,d^i_2,\cdots\}$ of Borel measurable functions $d^i_t:S_i\to A_i~(t\in T)$ %satisfying $\phi^i_t(A_i|x)=1(\forall x\in S)$ 
		such that $\pi^i_t(d^i_t(s_t^i)|{h}_t^i)= 1~(\forall t\in T,{h}_t^i=(s_1^i,a_1^i,s_2^i,a_2^i,\cdots,s_t^i)\in H_t^i)$, then $\pi^i$ is called a deterministic Markov strategy for player $i$.
	\end{enumerate}
\end{definition}

\begin{remark}
	The model considers incomplete information, since each player can only see their own state and action history. This complies with the aforementioned definition of strategy.
\end{remark}

The set of all (randomized) history-dependent strategies, all Markov strategies and all deterministic Markov strategies for player $i\in N$ are denoted by $\Pi_i$, $\Pi^i_M$ and $\Pi^i_{DM}$, respectively. And define $\boldsymbol{\Pi}:=\prod_{i\in N}\Pi_i, \boldsymbol{\Pi}_M:=\prod_{i\in N}\Pi^i_M,\boldsymbol{\Pi}_{DM}:=\prod_{i\in N}\Pi^i_{DM}$.

Based on the Ionescu-Tulcea theorem \cite[Appendix C.10]{hernandez1996}, for any initial system state $\boldsymbol{x}=(x^1,x^2,\cdots,x^n)\in \boldsymbol{S}$ and strategy profile $\boldsymbol{\pi}=(\pi^1,\pi^2,\cdots,\pi^n)\in \boldsymbol{\Pi}$ of all players, there exists a unique probability measure $\mathbb{P}_{\boldsymbol{x}}^{\boldsymbol{\pi}}$ (resp. 
$\mathbb{P}_{x^i}^{\pi^i}$ for each player $i$) on $(H_\infty,\mathscr{B}(H_\infty))$ (resp. 
$(H_\infty^i,\mathscr{B}(H_\infty^i))$) such that for any $\boldsymbol{s}=(s^1,s^2,\cdots,s^n)\in \boldsymbol{S},\boldsymbol{a}=(a^1,a^2,\cdots,a^n)\in \boldsymbol{A},i\in N$,
\begin{align}
	&\mathbb{P}_{\boldsymbol{x}}^{\boldsymbol{\pi}}\left(\tilde{\boldsymbol{X}}_1=\boldsymbol{x}\right)=1,\label{pm0}\\ 
	&\mathbb{P}_{\boldsymbol{x}}^{\boldsymbol{\pi}}\left(\tilde{\boldsymbol{A}}_t=\boldsymbol{a}\big|\tilde{\boldsymbol{X}}_1,\tilde{\boldsymbol{A}}_1,\cdots,\tilde{\boldsymbol{X}}_{t}\right)=%\prod_{i\in N}\pi_t^i\left(a^i\big|\tilde{X}_1,\tilde{A}_1,\cdots,\tilde{X}_{t}\right)=
	\prod_{i\in N}\pi_t^i\left(a^i\big|\tilde{X}_1^i,\tilde{A}_1^i,\cdots,\tilde{X}_{t}^i\right),\label{pm2}\\
	&\mathbb{P}_{\boldsymbol{x}}^{\boldsymbol{\pi}}\left(\tilde{\boldsymbol{X}}_{t+1}=\boldsymbol{s}\big|\tilde{\boldsymbol{X}}_1,\tilde{\boldsymbol{A}}_1,\cdots,\tilde{\boldsymbol{X}}_{t},\tilde{\boldsymbol{A}}_{t}\right)=q_t(\boldsymbol{s}|\tilde{\boldsymbol{X}}_{t},\tilde{\boldsymbol{A}}_{t}),\label{pm4}\\
	&\mathbb{P}_{x^i}^{\pi^i}\left(X_1^i=x^i\right)=1,\label{pm1}\\
	&\mathbb{P}_{x^i}^{\pi^i}\left(A_t^i=a^i|X_1^i,A_1^i,\cdots,X_{t}^i\right)=\pi_t^i\left(a^i|X_1^i,A_1^i,\cdots,X_{t}^i\right),\label{pm3}\\
	&\mathbb{P}_{x^i}^{\pi^i}\left(X_{t+1}^i=s^i\big|X_1^i,A_1^i,\cdots,X_{t}^i,A_{t}^i\right)=q_t^i(s^i|X_{t}^i,A_{t}^i),\label{pm5}
\end{align}
where $\{\tilde{\boldsymbol{X}}_t=(\tilde{X}_t^1,\tilde{X}_t^2,\cdots,\tilde{X}_t^n):t\in T\}$ (resp. 
$\{X_t^i:t\in T\}$) and $\{\tilde{\boldsymbol{A}}_t=(\tilde{A}_t^1,\tilde{A}_t^2,\cdots,\tilde{A}_t^n):t\in T\}$ (resp. 
$\{A_t^i:t\in T\}$) are respectively state and action processes. % Denote by $\mathbb{E}_x^{\{\pi}}$ the expectation operator with respect to $\mathbb{P}_x^{\{\pi}}$.
%\textbf{For convenience, we ignore the subscript of probability measure in the followings since the initial state is given.} 
Under (\ref{probabilityindependence}) and (\ref{pm0})-(\ref{pm5}), we can prove that
\begin{align}
	&\mathbb{P}^{\boldsymbol{\pi}}_{\boldsymbol{x}}\left(\tilde{\boldsymbol{X}}_t=\boldsymbol{s},\tilde{\boldsymbol{A}}_t=\boldsymbol{a}\right)\nonumber \\&=\mathbb{P}^{\boldsymbol{\pi}}_{\boldsymbol{x}}\left(\tilde{X}_t^1=s^1,\tilde{X}_t^2=s^2,\cdots,\tilde{X}_t^n=s^n,\tilde{A}_t^1=a^1,\tilde{A}_t^2=a^2,\cdots,\tilde{A}_t^n=a^n\right)\nonumber \\
	&=\prod_{i\in N}\mathbb{P}^{\pi^i}_{x^i}\left(X_t^i=s^i,A_t^i=a^i\right),\nonumber\\
	&~~~~~~~~~~~~~~~~~~~~~~~~~~~~\forall t\in T,\boldsymbol{s}=(s^1,s^2,\cdots,s^n)\in \boldsymbol{S},\boldsymbol{a}=(a^1,a^2,\cdots,a^n)\in \boldsymbol{A}.\nonumber
\end{align}

\begin{remark}
	It is possible to define the product probability measure of $\mathbb{P}^{\pi^i}_{x^i}$ directly, as shown in \cite{zhang2019,zhang2022}. However, under (\ref{probabilityindependence}), we provide the relations between the probability measure in centralized\footnote{Now, the centralized stochastic game represents that $\{\tilde{\boldsymbol{X}}_t:t\in T\}$ is the single (Markov) chain and all players jointly control it, when the strategy profile $\boldsymbol{\pi}\in \boldsymbol{\Pi}~(\boldsymbol{\pi}\in \boldsymbol{\Pi}_M)$ of all players is given.} and decentralized stochastic games here, i.e., $\mathbb{P}^{\boldsymbol{\pi}}_{\boldsymbol{x}}$ is (setwise) equal to $\prod_{i\in N}\mathbb{P}^{\pi^i}_{x^i}$.
\end{remark}

We consider the players with subjective behavior when facing randomness and uncertainty, i.e., prospect players. For the distortion of probability in PT, assume that each player has a subjective evaluation of other players' probabilities only. As their probability can be determined by their own strategy, whereas the strategies of others are uncertain. Thus, we have the ``probability measure" of player $i$ under PT \cite{etesami}:
\begin{align}
	&\mathbb{P}_{i,\rm{PT}}^{\boldsymbol{x},\boldsymbol{\pi}}\left(\tilde{\boldsymbol{X}}_t=\boldsymbol{s},\tilde{\boldsymbol{A}}_t=\boldsymbol{a}\right)\nonumber \\
	&:=\mathbb{P}^{\pi^i}_{x^i}\left(X_t^i=s^i,A_t^i=a^i\right)w_i\left(\prod_{j\in N,j\neq i}\mathbb{P}^{\pi^j}_{x^j}\left(X_t^j=s^j,A_t^j=a^j\right)\right),\nonumber
\end{align}
where $w_i:[0,1]\to [0,1]$ is a continuous weighting function for player $i$, e.g., $w_i(y)=e^{-(-\ln{y})^{\alpha}}(\alpha\in (0,1])$ \cite{prelec}.

Next, we define the expected prospect for player $i$.
\begin{definition}
	For each initial state $\boldsymbol{x}\in \boldsymbol{S}$ and strategy profile $\boldsymbol{\pi}\in \boldsymbol{\Pi}$, the expected prospect of player $i$ at time $t$ is 
	\begin{align}
		\mathbb{E}_{i,\rm{PT}}^{\boldsymbol{x},\boldsymbol{\pi}}\left[r_i(\tilde{\boldsymbol{X}}_t,\tilde{\boldsymbol{A}}_t)\right] := \sum_{\boldsymbol{s}\in \boldsymbol{S},\boldsymbol{a}\in \boldsymbol{A}}v_i\left(r_i(\boldsymbol{s},\boldsymbol{a})\right)\mathbb{P}_{i,\rm{PT}}^{\boldsymbol{x},\boldsymbol{\pi}}\left(\tilde{\boldsymbol{X}}_t=\boldsymbol{s},\tilde{\boldsymbol{A}}_t=\boldsymbol{a}\right),\nonumber
	\end{align}
	where $v_i:\mathbb{R}\to \mathbb{R}$ is the valuation function for player $i$, e.g., 
	$v_i(y)=\left \{
	\begin{array}{ll}
		y^{c_1},      &  y\geq 0, \\
		-c_2(-y)^{c_3},      & y<0,
	\end{array}
	\right.$
	($c_1,c_2,c_3>0$) \cite{al}.
\end{definition}

Now, we can introduce the optimality criterion: the expected discounted criterion under PT. For each initial state $\boldsymbol{x}=(x^1,x^2,\cdots,x^n)\in \boldsymbol{S}$ and strategy profile $\boldsymbol{\pi}=(\pi^1,\pi^2,\cdots,\pi^n)\in \boldsymbol{\Pi}$, the expected discounted payoff function under PT for player $i\in N$ is 
\begin{align}\label{criterion}
	V_i(\boldsymbol{x},\pi^i,\boldsymbol{\pi}^{-i}):=&\sum_{t=1}^{\infty}\beta^{t-1}\mathbb{E}_{i,\rm{PT}}^{\boldsymbol{x},\boldsymbol{\pi}}\left[r_i(\tilde{\boldsymbol{X}}_t,\tilde{\boldsymbol{A}}_t)\right]\nonumber\\
	=&\sum_{t=1}^{\infty}\beta^{t-1}\sum_{\boldsymbol{s}\in \boldsymbol{S},\boldsymbol{a}\in \boldsymbol{A}}v_i\left(r_i(\boldsymbol{s},\boldsymbol{a})\right)\mathbb{P}_{i,\rm{PT}}^{\boldsymbol{x},\boldsymbol{\pi}}\left(\tilde{\boldsymbol{X}}_t=\boldsymbol{s},\tilde{\boldsymbol{A}}_t=\boldsymbol{a}\right)\nonumber\\
	=&\sum_{t=1}^{\infty}\beta^{t-1}\sum_{\substack{s^1\in S_1,s^2\in S_2,\cdots,s^n\in S_n,\\a^1\in A_1,a^2\in A_2,\cdots,a^n\in A_n}}v_i\left(r_i(s^1,s^2,\cdots,s^n,a^1,a^2,\cdots,a^n)\right)\cdot \nonumber\\
	&\quad \mathbb{P}^{\pi^i}_{x^i}\left(X_t^i=s^i,A_t^i=a^i\right) w_i\left(\prod_{j\in N,j\neq i}\mathbb{P}^{\pi^j}_{x^j}\left(X_t^j=s^j,A_t^j=a^j\right)\right).
\end{align}
%in which $\mathbb{E}_{i,\rm{PT}}^{x,\pi}$ is the expectation operator with respect to $\mathbb{P}_{i,\rm{PT}}^{x,\pi}$, and $v_i$ is the valuation function for prosumer $i$.
Here, $(\gamma^i,{\boldsymbol{\pi}^{-i}})$ denotes the strategy profile that player $i$ takes strategy $\gamma^i$ and player $j~(\forall j\in N,j\neq i)$ takes strategy $\pi^j$. We give the definition of the optimization goal for stochastic game (\ref{model}) as below.

\begin{definition}
	Given an initial state $\boldsymbol{x}\in \boldsymbol{S}$, a strategy profile ${\boldsymbol{\pi}_*}=(\pi_*^1,\pi_*^2,\cdots,\pi_*^n)\in \boldsymbol{\Pi}$ is called a $\boldsymbol{x}$-Nash equilibrium if 
	\be \label{nashdef}
	V_i(\boldsymbol{x},\pi_*^i,\boldsymbol{\pi}_*^{-i}) \geq V_i(\boldsymbol{x},\pi^i,\boldsymbol{\pi}_*^{-i}), ~~\forall i\in N,\pi^i \in \Pi_i, 
	\ee
	i.e.,
	\be
	V_i(\boldsymbol{x},\pi_*^i,\boldsymbol{\pi}_*^{-i}) =\max_{\pi^i\in \Pi_i} V_i(\boldsymbol{x},\pi^i,\boldsymbol{\pi}_*^{-i}), ~~\forall i\in N. \nonumber
	\ee
	Moreover, if (\ref{nashdef}) holds for all initial state $\boldsymbol{x}\in \boldsymbol{S}$, then the strategy profile ${\boldsymbol{\pi}_*}$ is called a Nash equilibrium.
	%	For each $\varepsilon\geq 0$, a strategy profile ${\pi^*}=(\pi^*_1,\pi^*_2,\cdots,\pi^*_m)\in {\Pi}$ is called an $\varepsilon$-equilibrium if 
	%	\be
	%	 J(i,x,{\pi^*}) \geq J(i,x,[\pi_i,{\pi^*_{-i}}])-\varepsilon, ~~\forall i\in I,x\in S,\pi_i \in \Pi_i. \nonumber
	%	\ee
	%% Is the following definition too strict?
	%	\be
	%	J(i,x,[\pi_i,\pi^*_{-i}])-\varepsilon \leq J(i,x,\pi^*) \leq J(i,x,[\pi_i,\pi^*_{-i}])+\varepsilon, ~~\forall i\in I,x\in S,\pi_i \in \Pi_i. 
	%	\ee
\end{definition}

\section{The Existence of Nash Equilibria}\label{proofsection}

In this section, we prove the existence of Nash equilibria for stochastic game (\ref{model}).

Fix any $\boldsymbol{x}=(x^1,x^2,\cdots,x^n)\in \boldsymbol{S}$ and $\boldsymbol{\pi}=(\pi^1,\pi^2,\cdots,\pi^n)\in \boldsymbol{\Pi}$, we define that
\begin{align}
	\rho^i_t(s^i,a^i):=\mathbb{P}^{\pi^i}_{x^i}\left(X_t^i=s^i,A_t^i=a^i\right),\nonumber
\end{align}
for all $i\in N,t\in T,s^i\in S_i,a^i\in A_i$, which is the marginal distribution on the state-action pairs at any time. Let $\rho^i_t:=\left(\rho^i_t(s^i,a^i)\right)_{s^i\in S_i,a^i\in A_i}$ be a $|S_i|\times |A_i|$ matrix. Denote by $\rho^i:=(\rho^i_t:t\in T)$ the matrix vector. Thus, $\rho^i$ is determined by $x^i\in S_i$ and $\pi^i\in \Pi_i$. Now, note that the expected discounted payoff function under PT (\ref{criterion}) can be changed into
\begin{align}\label{defineJ}
	V_i(\boldsymbol{x},\pi^i,\boldsymbol{\pi}^{-i})
	&=\sum_{t=1}^{\infty}\beta^{t-1}\sum_{\substack{s^1\in S_1,s^2\in S_2,\cdots,s^n\in S_n,\\a^1\in A_1,a^2\in A_2,\cdots,a^n\in A_n}}v_i\left(r_i(s^1,s^2,\cdots,s^n,a^1,a^2,\cdots,a^n)\right)\cdot \nonumber\\
	&\quad\quad \mathbb{P}^{\pi^i}_{x^i}\left(X_t^i=s^i,A_t^i=a^i\right) w_i\left(\prod_{j\in N,j\neq i}\mathbb{P}^{\pi^j}_{x^j}\left(X_t^j=s^j,A_t^j=a^j\right)\right)\nonumber\\
	&=\sum_{t=1}^{\infty}\beta^{t-1}\sum_{s^i\in S_i,a^i\in A_i}\rho^i_t(s^i,a^i) \sum_{\substack{s^j\in S_j,a^j\in A_j,\\j\in N,j\neq i}}w_i\left(\prod_{j\in N,j\neq i}\rho^j_t(s^j,a^j)\right)\cdot\nonumber\\
	&\quad\quad v_i\left(r_i(s^1,s^2,\cdots,s^n,a^1,a^2,\cdots,a^n)\right)\nonumber\\
	&=:J_i(\rho^i,\boldsymbol{\rho}^{-i}),
\end{align}
in which $\boldsymbol{\rho}^{-i}$ is $(\rho^{j}:j\in N,j\neq i)$. Define that
\begin{align}\label{M}
	M_{\boldsymbol{x}}^i:=&\bigcup_{\pi^i\in \Pi_M^i}\bigg\{\rho^i=\left(\left(\rho^i_t(s^i,a^i)\right)_{s^i\in S_i,a^i\in A_i}:t\in T\right)\bigg|\rho^i_t(s^i,a^i)=\mathbb{P}^{\pi^i}_{x^i}\left(X_t^i=s^i,A_t^i=a^i\right),\nonumber\\
	&~~~~~~~~~~~~~~~~~~~~~~~~~~~~~~~~~~~~~~~~~~~~~~~~~~~~~~~~~~~~~~~~~~~~~~~~~\forall t\in T,s^i\in S_i,a^i\in A_i\bigg\}.
%mistake	M_i:=&\bigg\{\rho^i=(\rho^i_t:t\in T)\bigg|\rho^i_t(s^i,a^i)=\mathbb{P}^{\pi^i}\left(X_t^i=s^i,A_t^i=a^i\right),\nonumber\\
%	&~~~~~~~~~~~~~~~~~~~~~~~~~~~~~~~~~~~~~~~~~~~~~~~~~~~~~~~~~~~~~~~~~~~~~~\forall t\in T,s^i\in S_i,a^i\in A_i,\pi^i\in \Pi_i\bigg\}.\nonumber
\end{align}
It can be seen that $M_{\boldsymbol{x}}^i$ is nonempty, since the state and action spaces are nonempty.

When $\boldsymbol{x}\in \boldsymbol{S}$ and $\boldsymbol{\pi}^{-i}\in \prod_{j\in N,j\neq i}\Pi_j~(\forall i\in N)$ are fixed, $\boldsymbol{\rho}^{-i}$ is also fixed, and we have
\begin{align}\label{maxequal}
	\sup_{\rho^i\in M_{\boldsymbol{x}}^i}J_i(\rho^i,\boldsymbol{\rho}^{-i})=\sup_{\pi^i\in \Pi_M^i}V_i(\boldsymbol{x},\pi^i,\boldsymbol{\pi}^{-i})=\sup_{\pi^i\in \Pi_i}V_i(\boldsymbol{x},\pi^i,\boldsymbol{\pi}^{-i}),
\end{align} 
which the last equality is due to the existence of a Markov strategy with the same expected discounted payoff under PT for any given randomized history-dependent strategy\footnote{Fix initial state $\boldsymbol{x}=(x^1,x^2,\cdots,x^n)\in \boldsymbol{S}$, for any randomized history-dependent strategy $\boldsymbol{\pi}=(\pi^1,\pi^2,\cdots,\pi^n)\in \boldsymbol{\Pi}$, there exists a Markov strategy $\boldsymbol{\Psi}=(\Psi^1,\Psi^2,\cdots,\Psi^n)\in \boldsymbol{\Pi}_M$ such that
	\begin{align}
		\mathbb{P}^{\pi^i}_{x^i}\left(X_t^i=s^i,A_t^i=a^i\right)=\mathbb{P}^{\Psi^i}_{x^i}\left(X_t^i=s^i,A_t^i=a^i\right),~~\forall i\in N, t\in T,s^i\in S_i,a^i\in A_i,\nonumber
\end{align}which the proof is the nonstationary case of \cite[Theorem 5.5.1]{puterman}. Then, (\ref{criterion}) implies that the last equality of (\ref{maxequal}) holds.}.
If one of the supremums can be obtained in (\ref{maxequal}), then all the supremums can.%$\max_{\rho^i\in M_x^i}J_i(\rho^i,\rho^{-i})=\max_{\pi^i\in \Pi_M^i}V_i(x,\pi^i,{\pi}^{-i})=\max_{\pi^i\in \Pi_i}V_i(x,\pi^i,{\pi}^{-i})$.

For any $\boldsymbol{x}\in \boldsymbol{S}$, if we can verify that the static game\footnote{The players set is still $N$. The action space of player $i\in N$ is $M_{\boldsymbol{x}}^i$.} with payoffs $J_i(\rho^i,\boldsymbol{\rho}^{-i})$ has a Nash equilibrium $\boldsymbol{\rho}_{*,\boldsymbol{x}}=(\rho_{*,\boldsymbol{x}}^1,\rho_{*,\boldsymbol{x}}^2,\cdots,\rho_{*,\boldsymbol{x}}^n)\in\prod_{i\in N}M_{\boldsymbol{x}}^i$, stochastic game (\ref{model}) with payoffs $V_i(\boldsymbol{x},\pi^i,\boldsymbol{\pi}^{-i})$ also exists. Note that for any $\boldsymbol{x}\in \boldsymbol{S},i\in N$, we can at least find a $\pi_{*,\boldsymbol{x}}^i\in \Pi_M^i$ such that $\rho^i_{*,\boldsymbol{x},t}(s^i,a^i)=\mathbb{P}^{\pi^i_{*,\boldsymbol{x}}}_{x^i}\left(X_t^i=s^i,A_t^i=a^i\right)(\forall t\in T,s^i\in S_i,a^i\in A_i)$, then by (\ref{defineJ}) and (\ref{maxequal}), we have
\begin{align}\label{equivalent}
	V_i(\boldsymbol{x},\pi_{*,\boldsymbol{x}}^i,\boldsymbol{\pi}_{*,\boldsymbol{x}}^{-i})=J_i(\rho_{*,\boldsymbol{x}}^i,\boldsymbol{\rho}_{*,\boldsymbol{x}}^{-i})=\max_{\rho^i\in M_{\boldsymbol{x}}^i}J_i(\rho^i,\boldsymbol{\rho}_{*,\boldsymbol{x}}^{-i})=\max_{\pi^i\in \Pi_i}V_i(\boldsymbol{x},\pi^i,\boldsymbol{\pi}_{*,\boldsymbol{x}}^{-i}),
\end{align}
holds for all $i\in N$, i.e., $\boldsymbol{\pi}_{*,\boldsymbol{x}}=\left(\pi_{*,\boldsymbol{x}}^i,\boldsymbol{\pi}_{*,\boldsymbol{x}}^{-i}\right)\in \boldsymbol{\Pi}_M$ is a $\boldsymbol{x}$-Nash equilibrium of stochastic game (\ref{model}). Furthermore, we define a strategy $\boldsymbol{\pi}_*\in\boldsymbol{\Pi}$ by
\begin{align}\label{defNE}
	\boldsymbol{\pi}_*:=\sum_{\boldsymbol{x}\in \boldsymbol{S}}\mathbb{I}(\tilde{\boldsymbol{X}}_1=\boldsymbol{x})\boldsymbol{\pi}_{*,\boldsymbol{x}},
\end{align}
where $\mathbb{I}(\cdot)$ is the indicator function. Then, it is obviously that $\boldsymbol{\pi}_*$ is a Nash equilibrium of stochastic game (\ref{model}) with payoffs $V_i(\boldsymbol{x},\pi^i,\boldsymbol{\pi}^{-i})$.
Thus, we just need to prove that there exists a Nash equilibrium of the static game with payoffs $J_i(\rho^i,\boldsymbol{\rho}^{-i})$ for any fixed initial state. \textbf{In the followings, we assume that the initial state $\boldsymbol{x}\in \boldsymbol{S}$ is arbitrarily given.}

\begin{lemma}\label{lemma1}
	If the following statements hold for all $i\in N$, 
	\begin{enumerate}[(a)]
		\item $M_{\boldsymbol{x}}^i$ is a convex and closed set;
		\item Fix $\boldsymbol{\rho}^{-i}\in \prod_{j\in N,j\neq i}M_{\boldsymbol{x}}^j$, $J_i(\rho^i,\boldsymbol{\rho}^{-i})$ is concave in $\rho^i\in M_{\boldsymbol{x}}^i$;
		\item $J_i(\boldsymbol{\rho}):=J_i(\rho^i,\boldsymbol{\rho}^{-i})$ is continuous in $\boldsymbol{\rho}\in \prod_{i\in N}M_{\boldsymbol{x}}^i$,
	\end{enumerate}
    then there exists a Nash equilibrium of the static game with payoffs $J_i(\rho^i,\boldsymbol{\rho}^{-i})$.
\end{lemma}
{\it Proof:} Denote by $M_{\boldsymbol{x}}:=\prod_{i\in N}M_{\boldsymbol{x}}^i$. The correspondance $\Gamma^i:M_{\boldsymbol{x}}\twoheadrightarrow M_{\boldsymbol{x}}^i$ is defined by
\begin{align}
	\Gamma^i \boldsymbol{\rho}:=\left\{\eta^i\in M_{\boldsymbol{x}}^i\bigg|J_i(\eta^i,\boldsymbol{\rho}^{-i})=\max_{\xi^i\in M_{\boldsymbol{x}}^i}J_i(\xi^i,\boldsymbol{\rho}^{-i})\right\}.\nonumber
\end{align}
And define the correspondance $\Gamma:M_{\boldsymbol{x}}\twoheadrightarrow M_{\boldsymbol{x}}$ by $\Gamma \boldsymbol{\rho}:=\prod_{i\in N}\Gamma^i \boldsymbol{\rho}$. 

Note that the space of all real sequences $\prod_{n=1}^{\infty}\mathbb{R}$ is metrizable by \cite[Theorem 3.36]{aliprantis}, which the metric is defined by $d(\{y_n\},\{z_n\}):=\sum_{n=1}^{\infty}\frac{1}{2^n}\frac{|y_n-z_n|}{1+|y_n-z_n|}~(\forall \{y_n\},\{z_n\}\in \prod_{n=1}^{\infty}\mathbb{R})$. And $d$-convergence in $\prod_{n=1}^{\infty}\mathbb{R}$ is equivalent to pointwise convergence in $\mathbb{R}$. Thus, $M_{\boldsymbol{x}}\subset \prod_{n=1}^{\infty}\mathbb{R}$ and $\prod_{n=1}^{\infty}\mathbb{R}$ is a locally convex Hausdorff topological linear space. Denote by $E:=\{y=(y_1,y_2,\cdots)\in \prod_{n=1}^{\infty}\mathbb{R}|y_n\in [0,1],\forall n=1,2,\cdots\}$. Since $E$ is compact in the product topology (induced by metric $d$) by the Tychonoff theorem \cite[Theorem 2.61]{aliprantis}, $M_{\boldsymbol{x}}^i\subset E\subset \prod_{n=1}^{\infty}\mathbb{R}~(\forall i\in N)$ and (a) imply that $M_{\boldsymbol{x}}^i$ is nonempty, compact and convex for any $i\in N$, as well as $M_{\boldsymbol{x}}$. 

%mistake
%Let $l^{\infty}$ be the space of all bounded sequences (with metric induced by supremum norm), thus $M\subset l^{\infty}$ and $l^{\infty}$ is a locally convex Hausdorff topological linear space. Denote by $E:=\{x=(x_1,x_2,\cdots)\in l^{\infty}|x_i\in [0,1],\forall i=1,2,\cdots\}$. Since $E$ is compact (pointwise convergence) by the Tychonoff theorem, $M_i\subset E\subset l^{\infty}(\forall i\in N)$ and (a) imply that $M_i$ is nonempty, compact and convex for any $i\in N$, as well as $M$. 

For any $i\in N,\boldsymbol{\rho}\in M_{\boldsymbol{x}}$, since $M_{\boldsymbol{x}}^i$ is compact, (c) implies that there exists $\zeta^i\in M_{\boldsymbol{x}}^i$ such that $J_i(\zeta^i,\boldsymbol{\rho}^{-i})=\max_{\xi^i\in M_{\boldsymbol{x}}^i}J_i(\xi^i,\boldsymbol{\rho}^{-i})$, thus $\Gamma^i$ has nonempty values.
For each $\eta^i,\zeta^i\in \Gamma^i \boldsymbol{\rho}$, and $\lambda\in [0,1]$, we have $\lambda \eta^i+(1-\lambda)\zeta^i\in M_{\boldsymbol{x}}^i$ by the convexity of $M_{\boldsymbol{x}}^i$. And (b) implies that
\begin{align}
	J_i(\lambda \eta^i+(1-\lambda)\zeta^i,\boldsymbol{\rho}^{-i})&\geq \lambda J_i(\eta^i,\boldsymbol{\rho}^{-i})+(1-\lambda)J_i(\zeta^i,\boldsymbol{\rho}^{-i})\nonumber\\
	&=\lambda\max_{\xi^i\in M_{\boldsymbol{x}}^i}J_i(\xi^i,\boldsymbol{\rho}^{-i})+(1-\lambda)\max_{\xi^i\in M_{\boldsymbol{x}}^i}J_i(\xi^i,\boldsymbol{\rho}^{-i})\nonumber\\
	&=\max_{\xi^i\in M_{\boldsymbol{x}}^i}J_i(\xi^i,\boldsymbol{\rho}^{-i}),\nonumber
\end{align}
then $\lambda \eta^i+(1-\lambda)\zeta^i\in \Gamma^i \boldsymbol{\rho}$, i.e., $\Gamma^i$ has convex values. Thus, $\Gamma$ has nonempty convex values.

Denote the graph of the correspondence $\Gamma$ by ${\rm Gr}\Gamma:=\{(\boldsymbol{\rho},\boldsymbol{\eta})|\boldsymbol{\rho} \in M_{\boldsymbol{x}},\boldsymbol{\eta}\in\Gamma\boldsymbol{\rho}\}$. For any sequences $\{(\boldsymbol{\rho}_k,\boldsymbol{\eta}_k),k\geq 1\}\subset {\rm Gr}\Gamma$, if $(\boldsymbol{\rho}_k,\boldsymbol{\eta}_k)\to (\boldsymbol{\rho},\boldsymbol{\eta})~(k\to \infty)$, we claim that $(\boldsymbol{\rho},\boldsymbol{\eta})\in {\rm Gr}\Gamma$. Note that $M_{\boldsymbol{x}}$ is compact and $\boldsymbol{\rho}_k \to \boldsymbol{\rho}$, thus $\boldsymbol{\rho} \in M_{\boldsymbol{x}}$. Moreover, $\boldsymbol{\eta}_k\in \Gamma \boldsymbol{\rho}_k~(k\geq 1)$, that is, for any $i\in N,\xi^i\in M_{\boldsymbol{x}}^i$, $J_i(\eta_k^i,\boldsymbol{\rho}_k^{-i})\geq J_i(\xi^i,\boldsymbol{\rho}_k^{-i})$. Let $k\to\infty$, we have $J_i(\eta^i,\boldsymbol{\rho}^{-i})\geq J_i(\xi^i,\boldsymbol{\rho}^{-i})$ by (c). Since $i$ and $\xi^i$ are arbitrary, we verify that $\boldsymbol{\eta}\in \Gamma\boldsymbol{\rho}$. Then, $\Gamma$ has closed graph.

Now, by the Kakutani-Fan-Glicksberg fixed point theorem \cite[Corollary 17.55]{aliprantis}, there exists $\boldsymbol{\rho}_*\in M_{\boldsymbol{x}}$ such that $\boldsymbol{\rho}_*\in\Gamma \boldsymbol{\rho}_*$, i.e.,
\begin{align}%\label{fixpoint}
	J_i(\rho_*^i,\boldsymbol{\rho}_*^{-i})\geq J_i(\xi^i,\boldsymbol{\rho}_*^{-i}),~~\forall i\in N,\xi^i\in M_{\boldsymbol{x}}^i.\nonumber
\end{align} 
Thus, $\boldsymbol{\rho}_*$ is a Nash equilibrium of the static game with payoffs $J_i(\rho^i,\boldsymbol{\rho}^{-i})$.

\te

Before verifying the conditions of Lemma \ref{lemma1}, we give an equivalent statement of $M_{\boldsymbol{x}}^i$ for any $i\in N$.

\begin{lemma}\label{lemma2}
	$\rho^i\in M_{\boldsymbol{x}}^i$ if and only if $\rho^i = \left(\rho^i_t(s^i,a^i)\right)_{s^i\in S_i,a^i\in A_i,t\in T}$ is an nonegative sequence satisfying 
	\begin{align}\label{lemma2.1}
		\sum_{a^i\in A_i} \rho^i_1(s^i,a^i) = \delta_{x^i}(s^i),
	\end{align}
and
    \begin{align}\label{lemma2.2}
    	\sum_{a^i\in A_i} \rho^i_{t+1}(s^i,a^i) =\sum_{s^i_t\in S_i, a^i_t\in A_i} \rho^i_t(s^i_t,a^i_t) q^i_t(s^i|s^i_t,a^i_t),~~t\in T,
    \end{align}
which $\delta_{x^i}$ is the 
Dirac measure concentrated at $x^i$.
\end{lemma}
{\it Proof:} If $\rho^i\in M_{\boldsymbol{x}}^i$, there exists $\pi^i\in\Pi_M^i$ such that 
$\rho^i_t(s^i,a^i)=\mathbb{P}^{\pi^i}_{x^i}\left(X_t^i=s^i,A_t^i=a^i\right)\in [0,1]$ for all $t\in T,s^i\in S_i,a^i\in A_i$. Then, 
\begin{align}
	\sum_{a^i\in A_i} \rho^i_1(s^i,a^i) = \sum_{a^i\in A_i}\mathbb{P}^{\pi^i}_{x^i}\left(X_1^i=s^i,A_1^i=a^i\right) = \mathbb{P}^{\pi^i}_{x^i}\left(X_1^i=s^i\right) =  \delta_{x^i}(s^i),\nonumber
\end{align}
and 
\begin{align}
	\sum_{a^i\in A_i} \rho^i_{t+1}(s^i,a^i) &= \sum_{a^i\in A_i} \mathbb{P}^{\pi^i}_{x^i}\left(X_{t+1}^i=s^i,A_{t+1}^i=a^i\right) \nonumber\\
	&= \mathbb{P}^{\pi^i}_{x^i}\left(X_{t+1}^i=s^i\right)\nonumber\\ 
	&= \sum_{s^i_t\in S_i, a^i_t\in A_i} \mathbb{P}^{\pi^i}_{x^i}\left(X_t^i=s_t^i,A_t^i=a_t^i\right) \mathbb{P}^{\pi^i}_{x^i}\left(X_{t+1}^i=s^i \big|X_t^i=s_t^i,A_t^i=a_t^i\right) \nonumber\\
	&=\sum_{s^i_t\in S_i, a^i_t\in A_i} \rho^i_t(s^i_t,a^i_t) q^i_t(s^i|s^i_t,a^i_t),~~t\in T.\nonumber
\end{align}

If nonegative sequence $\rho^i = \left(\rho^i_t(s^i,a^i)\right)_{s^i\in S_i,a^i\in A_i,t\in T}$ satisfies (\ref{lemma2.1}) and (\ref{lemma2.2}), we need to find $\pi^i\in\Pi_M^i$ such that 
\begin{align}\label{lemma2.proof}
	\rho^i_t(s^i,a^i)=\mathbb{P}^{\pi^i}_{x^i}\left(X_t^i=s^i,A_t^i=a^i\right),~~\forall t\in T,s^i\in S_i,a^i\in A_i.
\end{align}
Define that
\begin{align}\label{rhotopi}
	\pi_t^i(a^i|s^i) = \left \{
	\begin{array}{ll}
		\frac{\rho^i_t(s^i,a^i)}{\sum_{a^i\in A_i}\rho^i_t(s^i,a^i)},      &  \sum_{a^i\in A_i}\rho^i_t(s^i,a^i)>0, \\
		\frac{1}{|A_i|},      & \sum_{a^i\in A_i}\rho^i_t(s^i,a^i)=0.
	\end{array}
	\right.
\end{align} 
For every $t\in T,s^i\in S_i$, it can be seen that $\sum_{a^i\in A_i} \pi_t^i(a^i|s^i) =1$, and $\pi_t^i(\cdot|s^i)\in \mathscr{P}(A_i)$, %finite state space, 
thus, we have $\pi^i\in\Pi_M^i$. In the followings, we prove that $\pi^i$ defined in (\ref{rhotopi}) satisfies (\ref{lemma2.proof}) by induction. 

First, by (\ref{lemma2.1}) and (\ref{rhotopi}), 
\begin{align}
	\mathbb{P}^{\pi^i}_{x^i}\left(X_1^i=s^i,A_1^i=a^i\right) &= \delta_{x^i}(s^i) \pi_1^i(a^i|s^i) \nonumber\\
	&=\pi_1^i(a^i|s^i) \sum_{a^i\in A_i} \rho^i_1(s^i,a^i)\nonumber\\
	&=\left \{
	\begin{array}{ll}
		\rho^i_1(s^i,a^i),&\sum_{a^i\in A_i} \rho^i_1(s^i,a^i)>0,\\
		0=\rho^i_1(s^i,a^i),&\sum_{a^i\in A_i} \rho^i_1(s^i,a^i)=0.
	\end{array}\right.\nonumber
\end{align}
Assume that (\ref{lemma2.proof}) holds at time $t$, then (\ref{lemma2.2}) and (\ref{rhotopi}) imply that 
\begin{align}
	\mathbb{P}^{\pi^i}_{x^i}\left(X_{t+1}^i=s^i,A_{t+1}^i=a^i\right) &=\sum_{s^i_t\in S_i, a^i_t\in A_i}\mathbb{P}^{\pi^i}_{x^i}\left(X_t^i=s_t^i,A_t^i=a_t^i\right)\cdot\nonumber\\ &~~~~~~~~~~~~~~~~~~~~~~~\mathbb{P}^{\pi^i}_{x^i}\left(X_{t+1}^i=s^i|X_t^i=s_t^i,A_t^i=a_t^i\right)\cdot\nonumber\\ &~~~~~~~~~~~~~~~~~~~~~~~\mathbb{P}^{\pi^i}_{x^i}\left(A_{t+1}^i=a^i|X_t^i=s_t^i,A_t^i=a_t^i,X_{t+1}^i=s^i\right)\nonumber\\
	&=\sum_{s^i_t\in S_i, a^i_t\in A_i} \rho^i_t(s_t^i,a_t^i) q_t^i(s^i|s^i_t, a^i_t) \pi_{t+1}^i(a^i|s^i)\nonumber\\
	&=\pi_{t+1}^i(a^i|s^i) \sum_{a^i\in A_i}\rho^i_{t+1}(s^i,a^i)\nonumber\\
	&=\left \{
	\begin{array}{ll}
		\rho^i_{t+1}(s^i,a^i),&\sum_{a^i\in A_i} \rho^i_{t+1}(s^i,a^i)>0,\\
		0=\rho^i_{t+1}(s^i,a^i),&\sum_{a^i\in A_i} \rho^i_{t+1}(s^i,a^i)=0.
	\end{array}\right.\nonumber
\end{align}
\te

\begin{remark}
	(\ref{rhotopi}) provides a construction method of strategy satisfying (\ref{lemma2.proof}) from $\rho^i\in M_{\boldsymbol{x}}^i$. Thus, we just need to consider the static nonzero-sum game with payoffs $J_i(\rho^i,\boldsymbol{\rho}^{-i})$ theoretically, and solve it by nonlinear programming \cite[Theorem 3.4.1]{barron}.
\end{remark}

\begin{theorem}\label{theorem1}
	There exists a Nash equilibrium for the nonstationary decentralized stochastic game (\ref{model}) with the expected discounted payoff function under PT. Especially, for a given initial state $\boldsymbol{x}\in \boldsymbol{S}$, there exists a Markov $\boldsymbol{x}$-Nash equilibrium.
\end{theorem}
{\it Proof:} By Lemma \ref{lemma1}, (\ref{equivalent}) and (\ref{defNE}), we just need to show that Lemma \ref{lemma1} (a)-(c) holds for each $i\in N$.

For any sequences $\{\rho^{k,i},k\geq 1\}\subset M_{\boldsymbol{x}}^i$ %(ignore superscript $i$ for convenience)
such that $\rho^{k,i}\to\rho^i ~(k\to\infty)$ (pointwise convergence). It is obvious that $\rho^i = \left(\rho^i_t(s^i,a^i)\right)_{s^i\in S_i,a^i\in A_i,t\in T}$ is a sequence on $[0,1]$.
By Lemma \ref{lemma2}, finite state and action spaces, we have
\begin{align}
	\sum_{a^i\in A_i} \rho^i_1(s^i,a^i) = \sum_{a^i\in A_i} \lim_{k\to \infty} \rho^{k,i}_1(s^i,a^i) = \lim_{k\to \infty} \sum_{a^i\in A_i} \rho^{k,i}_1(s^i,a^i) = \lim_{k\to \infty} \delta_{x^i}(s^i) = \delta_{x^i}(s^i),\nonumber
\end{align}
and
\begin{align}
	\sum_{a^i\in A_i} \rho^i_{t+1}(s^i,a^i) &= \lim_{k\to \infty} \sum_{a^i\in A_i} \rho^{k,i}_{t+1}(s^i,a^i) \nonumber\\
	&= \lim_{k\to \infty} \sum_{s^i_t\in S_i, a^i_t\in A_i} \rho^{k,i}_t(s^i_t,a^i_t) q^i_t(s^i|s^i_t,a^i_t) \nonumber\\
	&= \sum_{s^i_t\in S_i, a^i_t\in A_i} \lim_{k\to \infty} \rho^{k,i}_t(s^i_t,a^i_t) q^i_t(s^i|s^i_t,a^i_t) \nonumber\\
	&=  \sum_{s^i_t\in S_i, a^i_t\in A_i} \rho^i_t(s^i_t,a^i_t) q^i_t(s^i|s^i_t,a^i_t),~~t\in T.\nonumber
\end{align}
Thus, $\rho^i\in M_{\boldsymbol{x}}^i$ and $M_{\boldsymbol{x}}^i$ is closed.

For any $\rho^i,\eta^i\in M_{\boldsymbol{x}}^i$ and $\lambda\in [0,1]$, denote $\lambda\rho^i+(1-\lambda)\eta^i$ by $\xi^i$, then $\xi^i$ is still a sequence on $[0,1]$. By Lemma \ref{lemma2},
\begin{align}
	\sum_{a^i\in A_i} \xi^i_1(s^i,a^i) &= \sum_{a^i\in A_i} \left(\lambda \rho^i_1(s^i,a^i) + (1-\lambda)\eta^i_1(s^i,a^i) \right)\nonumber\\
	&=\lambda \sum_{a^i\in A_i} \rho^i_1(s^i,a^i) + (1-\lambda) \sum_{a^i\in A_i} \eta^i_1(s^i,a^i)\nonumber\\
	&= \lambda \delta_{x^i}(s^i) + (1-\lambda) \delta_{x^i}(s^i)\nonumber\\
	&=\delta_{x^i}(s^i),\nonumber
\end{align}
and 
\begin{align}
	\sum_{a^i\in A_i} \xi^i_{t+1}(s^i,a^i) &= \sum_{a^i\in A_i} \left(\lambda \rho^i_{t+1}(s^i,a^i) + (1-\lambda)\eta^i_{t+1}(s^i,a^i) \right)\nonumber\\
	&=\sum_{s^i_t\in S_i, a^i_t\in A_i}\lambda \rho^i_t(s^i_t,a^i_t) q^i_t(s^i|s^i_t,a^i_t) + \sum_{s^i_t\in S_i, a^i_t\in A_i}(1-\lambda) \eta^i_t(s^i_t,a^i_t) q^i_t(s^i|s^i_t,a^i_t)\nonumber\\
	&=\sum_{s^i_t\in S_i, a^i_t\in A_i} \left( \lambda \rho^i_t(s^i_t,a^i_t) + (1-\lambda) \eta^i_t(s^i_t,a^i_t)\right) q^i_t(s^i|s^i_t,a^i_t)\nonumber\\
	&=
	\sum_{s^i_t\in S_i, a^i_t\in A_i} \xi^i_t(s^i_t,a^i_t) q^i_t(s^i|s^i_t,a^i_t),~~t\in T.\nonumber
\end{align}
Then, $\xi^i\in M_{\boldsymbol{x}}^i$, i.e., $M_{\boldsymbol{x}}^i$ is convex.

Fix $\boldsymbol{\rho}^{-i}\in \prod_{j\in N,j\neq i}M_{\boldsymbol{x}}^j$, (\ref{defineJ}) implies that $J_i(\rho^i,\boldsymbol{\rho}^{-i})$ is a linear function with respect to $\rho^i\in M_{\boldsymbol{x}}^i$, thus it is a concave function. 

By the continuity of $w_i$ and the dominated convergence theorem, $J_i(\boldsymbol{\rho})=J_i(\rho^i,\boldsymbol{\rho}^{-i})$ is continuous in $\boldsymbol{\rho}\in \prod_{i\in N}M_{\boldsymbol{x}}^i$.
\te

\begin{remark}\label{thm_rem}
%	If $\pi^i,\psi^i\in\Pi_M^i$, we cannot construct a Markov strategy $\varphi^i$ in the form of (\ref{defstrategy}). That's why we consider the (randomized) history-dependent strategy in the definition of $M_x^i$. 
\begin{enumerate}[(a)]
	\item If $M_{\boldsymbol{x}}^i$ (\ref{M}) is defined by the union of all randomized history-dependent strategies $\pi^i\in \Pi_i$, the conditions of Lemma \ref{lemma1} still hold. Now, for any given initial state $\boldsymbol{x}\in \boldsymbol{S}$, there exists a randomized history-dependent $\boldsymbol{x}$-Nash equilibrium.
	\item The above method of solving stochastic game by static game is suitable for the finite horizon, i.e., time set $T:=\{1,2,\cdots,\bar{T}\}~(\bar{T}<\infty)$, and the expected finite horizon payoff function under PT for player $i$ is $\bar{V}_i(\boldsymbol{x},\pi^i,\boldsymbol{\pi}^{-i}):=\sum_{t=1}^{\bar{T}}\mathbb{E}_{i,\rm{PT}}^{\boldsymbol{x},\boldsymbol{\pi}}\left[r_i(\tilde{\boldsymbol{X}}_t,\tilde{\boldsymbol{A}}_t)\right]$.
	\item Theorem \ref{theorem1} can be obtained when the payoff function is nonstationary and uniformly bounded with respect to time $t$ (denoted by $r_t^i$) or when the action set for player $i$ depends on time $t$ and own state $s^i$ (denoted by $A_t^i(s^i)$). Let the weighting function $w_i(y)=y$ and the valuation function $v_i(y)=y$, then Theorem \ref{theorem1} establishes the existence of Nash equilibrium for the nonstationary decentralized stochastic games with the general expected discounted criterion (without PT) simultaneously. 
\end{enumerate}	
\end{remark}

\section{Algorithm}\label{algsection}

In this section, we still fix the initial state $\boldsymbol{x}=(x^1,x^2,\cdots,x^n)\in \boldsymbol{S}$, since it is easy to construct a Nash equilibrium from $\boldsymbol{x}$-Nash equilibrium by (\ref{defNE}). From Section \ref{proofsection}, we derive the existence of Markov $\boldsymbol{x}$-Nash equilibrium (denoted by $\boldsymbol{\Psi}_*$) of stochastic game (\ref{model}). However, the above proof of $\boldsymbol{\Psi}_*$ is not constructive based on the Kakutani-Fan-Glicksberg fixed point theorem. Considering the actual situation, it is difficult to find an exact Markov $\boldsymbol{x}$-Nash equilibrium since the set of Markov strategies has infinite elements. %it is not necessarily possible to find the Markov Nash equilibrium since the set of Markov strategies has infinite elements. 
Thus, this section aims to provide an algorithm to find a Markov $\varepsilon$-equilibrium. First, we introduce the definition of $\varepsilon$-equilibrium.

\begin{definition}
	For any $\varepsilon\geq 0$, a strategy profile ${\boldsymbol{\pi}_*}=(\pi_*^1,\pi_*^2,\cdots,\pi_*^n)\in \boldsymbol{\Pi}$ is called an $\varepsilon$-equilibrium if 
	\be
	V_i(\boldsymbol{x},\pi_*^i,\boldsymbol{\pi}_*^{-i}) \geq V_i(\boldsymbol{x},\pi^i,\boldsymbol{\pi}_*^{-i})-\varepsilon, ~~\forall i\in N,\pi^i \in \Pi_i. \nonumber
	\ee
	It is obviously that $\boldsymbol{\pi}_*$ is a $\boldsymbol{x}$-Nash equilibrium when $\varepsilon=0$.
\end{definition}

Before finding an $\varepsilon$-equilibrium ($\varepsilon>0$), we analyze the properties of the criterion function $V_i(\boldsymbol{x},\boldsymbol{\pi})=V_i(\boldsymbol{x},\pi^i,\boldsymbol{\pi}^{-i})$.
%Now, we analyze the properties of the criterion function $V_i(\pi)=V_i(\pi^i,{\pi}^{-i})$ to provide a basis for finding an $\varepsilon$-equilibrium ($\varepsilon>0$).

\begin{proposition}\label{unicont}
	For all $i\in N$, $V_i(\boldsymbol{x},\boldsymbol{\pi})$ is uniformly continuous on $\boldsymbol{\pi}=(\pi^1,\pi^2,\cdots,\pi^n)\in \boldsymbol{\Pi}_M$.
\end{proposition}
{\it Proof:} Note that $\boldsymbol{\Pi}_M=\prod_{i\in N}\prod_{t=1}^{\infty}\mathscr{P}(A_i)^{|S_i|}$ is a compact metric space endowed with product weak topology. Here, weak convergence of probability measures is equivalent to pointwise convergence due to the finite action space. Thus, we just need to verify that $V_i(\boldsymbol{x},\boldsymbol{\pi})$ is continuous on $\boldsymbol{\pi}=(\pi^1,\pi^2,\cdots,\pi^n)\in \boldsymbol{\Pi}_M$ for all $i\in N$.

For any sequences $\{\boldsymbol{\pi}^k,k\geq 1\}\subset \boldsymbol{\Pi}_M$ such that $\boldsymbol{\pi}^k\to \boldsymbol{\pi}\in \boldsymbol{\Pi}_M~(k\to\infty)$, i.e., $\pi^{k,i}_t(a^i|s^i)\xrightarrow{k\to\infty} \pi^i_t(a^i|s^i)$ holds for all $i\in N,t\in T,s^i\in S_i,a^i\in A_i$. Then,
\begin{align}\label{limit}
	\mathbb{P}^{\pi^{k,i}}_{x^i}\left(X_t^i=s^i,A_t^i=a^i\right)&=\sum_{\substack{s^i_2,s^i_3,\cdots,s^i_{t-1}\in S_i,\\a^i_1,a^i_2,\cdots,a^i_{t-1}\in A_i}}\pi^{k,i}_1(a^i_1|x^i)q_1^i(s^i_2|x^i,a^i_1)\pi^{k,i}_2(a^i_2|s^i_2)\cdots\nonumber\\ &~~~~~~~~~~~~~~~~~~~~~~~~~~~~~~q_{t-1}^i(s^i|s^i_{t-1},a^i_{t-1})\pi^{k,i}_t(a^i|s^i)\nonumber\\
	&\xrightarrow{k\to\infty}\sum_{\substack{s^i_2,s^i_3,\cdots,s^i_{t-1}\in S_i,\\a^i_1,a^i_2,\cdots,a^i_{t-1}\in A_i}}\pi^{i}_1(a^i_1|x^i)q_1^i(s^i_2|x^i,a^i_1)\pi^{i}_2(a^i_2|s^i_2)\cdots\nonumber\\ &~~~~~~~~~~~~~~~~~~~~~~~~~~~~~~~~~~~~~q_{t-1}^i(s^i|s^i_{t-1},a^i_{t-1})\pi^{i}_t(a^i|s^i)\nonumber\\
	&=\mathbb{P}^{\pi^{i}}_{x^i}\left(X_t^i=s^i,A_t^i=a^i\right).
\end{align}
Recall that the definition of $V_i(\boldsymbol{x},\boldsymbol{\pi})$ in (\ref{criterion}), (\ref{limit}), the continuity of $w_i$ and the dominated convergence theorem imply that $V_i(\boldsymbol{x},\boldsymbol{\pi}^k)\xrightarrow{k\to\infty}V_i(\boldsymbol{x},\boldsymbol{\pi})$.
\te

\begin{remark}\label{uni_con_his}
	Actually, the uniform continuity of $V_i(\boldsymbol{x},\boldsymbol{\pi})$ on $\boldsymbol{\pi}\in \boldsymbol{\Pi}$ can be similarly proved to Proposition \ref{unicont}.
\end{remark}

By proposition \ref{unicont}, if two Markov strategies $\boldsymbol{\pi},\boldsymbol{\psi}\in\boldsymbol{\Pi}_M$ are closed enough, for example, there exists $\delta=\delta(\varepsilon)>0$ such that 
\begin{align}\label{leqdelta}
	\sup_{i\in N,t\in T,s^i\in S_i,a^i\in A_i}|\pi^i_t(a^i|s^i)-\psi^i_t(a^i|s^i)|\leq \delta,
\end{align}
we have 
\begin{align}\label{aim}
|V_i(\boldsymbol{x},\boldsymbol{\pi})-V_i(\boldsymbol{x},\boldsymbol{\psi})|\leq \frac{\varepsilon}{2},~~\forall i\in N. 
\end{align}
For any $\boldsymbol{\Psi}\in \boldsymbol{\Pi}_M$ such that $\sup_{i\in N,t\in T,s^i\in S_i,a^i\in A_i}|{\Psi}^i_t(a^i|s^i)-\Psi^i_{*,t}(a^i|s^i)|\leq \delta$, we have
\begin{align}
	V_i(\boldsymbol{x},{\Psi}^i,\boldsymbol{\Psi}^{-i})+\frac{\varepsilon}{2}\geq V_i(\boldsymbol{x},\Psi_*^i,\boldsymbol{\Psi}_*^{-i})\geq V_i(\boldsymbol{x},\pi^i,\boldsymbol{\Psi}_*^{-i})\geq V_i(\boldsymbol{x},\pi^i,\boldsymbol{\Psi}^{-i})-\frac{\varepsilon}{2},\nonumber
\end{align}
holds for all $i\in N,\pi^i\in \Pi_i$, which the third inequality is from Remark \ref{uni_con_his}, then, $\boldsymbol{\Psi}$ is a Markov $\varepsilon$-equilibrium. Thus, we can divide the set of Markov strategies into finite strategies such that it is closed enough between strategies under the sense of (\ref{leqdelta}), up to find an $\varepsilon$-equilibrium under a sampling without replacement method. In the followings, we provide a sufficient condition of $\varepsilon$-equilibrium. First, we give a lemma.%In order to investigate if there is an $\varepsilon$-equilibrium in these finite strategies, we provide a sufficient condition in the followings. First, we give a lemma.%Furthermore, the followings is to verify if these finite strategies are $\varepsilon$-equilibrium. For this purpose, we propose a lemma.

\begin{lemma}\label{finiteappro}
	For any $\varepsilon>0$, there exists a constant $\tilde{m}=\tilde{m}(\varepsilon)>0$ such that
	\begin{align}
		\left|V_i(\boldsymbol{x},\boldsymbol{\pi})-V_i^{m}(\boldsymbol{x},\boldsymbol{\pi})\right|\leq \frac{\varepsilon}{3},~~\forall i\in N, \boldsymbol{\pi}\in\boldsymbol{\Pi},\nonumber
	\end{align}
    holds when $m\geq \tilde{m}$, where $V_i^{m}(\boldsymbol{x},\boldsymbol{\pi}):=\sum_{t=1}^{m}\beta^{t-1}\mathbb{E}_{i,\rm{PT}}^{\boldsymbol{x},\boldsymbol{\pi}}\left[r_i(\tilde{\boldsymbol{X}}_t,\tilde{\boldsymbol{A}}_t)\right]$ represents the $m$-step finite truncation of $V_i(\boldsymbol{x},\boldsymbol{\pi})$.
\end{lemma}
{\it Proof:} Let $K_i:=\sup_{\boldsymbol{s}\in \boldsymbol{S},\boldsymbol{a}\in \boldsymbol{A}}\left|v_i\left(r_i(\boldsymbol{s},\boldsymbol{a})\right)\right|$, then $K_i<\infty$ for all $i\in N$, since both state space and action space are finite. Now,
\begin{align}
	\left|V_i(\boldsymbol{x},\boldsymbol{\pi})-V_i^{m}(\boldsymbol{x},\boldsymbol{\pi})\right|&=\left|\sum_{t=m+1}^{\infty}\beta^{t-1}\mathbb{E}_{i,\rm{PT}}^{\boldsymbol{x},\boldsymbol{\pi}}\left[r_i(\tilde{\boldsymbol{X}}_t,\tilde{\boldsymbol{A}}_t)\right]\right|\nonumber\\
	&=\left|\sum_{t=m+1}^{\infty}\beta^{t-1}\sum_{\boldsymbol{s}\in \boldsymbol{S},\boldsymbol{a}\in \boldsymbol{A}}v_i\left(r_i(\boldsymbol{s},\boldsymbol{a})\right)\mathbb{P}_{i,\rm{PT}}^{\boldsymbol{x},\boldsymbol{\pi}}\left(\tilde{\boldsymbol{X}}_t=\boldsymbol{s},\tilde{\boldsymbol{A}}_t=\boldsymbol{a}\right)\right|\nonumber\\
	&\leq \sum_{t=m+1}^{\infty}\beta^{t-1}\sum_{\boldsymbol{s}\in \boldsymbol{S},\boldsymbol{a}\in \boldsymbol{A}}\left|v_i\left(r_i(\boldsymbol{s},\boldsymbol{a})\right)\right|\nonumber\\
	&\leq K_i |\boldsymbol{S}| |\boldsymbol{A}|\frac{\beta^m}{1-\beta},~~\forall i\in N,\boldsymbol{\pi}\in\boldsymbol{\Pi}.\nonumber
\end{align}
Thus, when $m\geq \left\lceil \max_{i\in N}\left\{\frac{\ln{\frac{(1-\beta)\varepsilon}{3K_i |\boldsymbol{S}| |\boldsymbol{A}|}}}{\ln{\beta}}\right\} \right\rceil =:\tilde{m}(\varepsilon)\footnote{$\lceil \cdot \rceil$ is the integer ceiling function.}$, we have $\left|V_i(\boldsymbol{x},\boldsymbol{\pi})-V_i^{m}(\boldsymbol{x},\boldsymbol{\pi})\right|\leq \frac{\varepsilon}{3}$ for all $i\in N,\boldsymbol{\pi}\in\boldsymbol{\Pi}$.
\te

From Lemma \ref{finiteappro}, if $\boldsymbol{\Psi}\in\boldsymbol{\Pi}_M$ satisfies
\begin{align}\label{finiteclose}
	V_i^{\tilde{m}}(\boldsymbol{x},\Psi^i,\boldsymbol{\Psi}^{-i})\geq V_i^{\tilde{m}}(\boldsymbol{x},\pi^i,\boldsymbol{\Psi}^{-i})-\frac{\varepsilon}{3},~~\forall i\in N,\pi^i\in\Pi_i,
\end{align}
then $\boldsymbol{\Psi}$ is a Markov $\varepsilon$-equilibrium for stochastic game (\ref{model}). In other words, (\ref{finiteclose}) is a sufficient condition of Markov $\varepsilon$-equilibrium.%to verify whether the Markov strategies are $\varepsilon$-equilibrium or not, we just need to consider (\ref{finiteclose}).

Fixed others strategies $\boldsymbol{\Psi}^{-i}\in \prod_{j\in N,j\neq i}\Pi_M^j$, note that 
\begin{align}
	&V_i^{\tilde{m}}(\boldsymbol{x},\pi^i,\boldsymbol{\Psi}^{-i})\nonumber\\
	&=\sum_{t=1}^{\tilde{m}}\sum_{s^i\in S_i,a^i\in A_i}\mathbb{P}^{\pi^i}_{x^i}\left(X_t^i=s^i,A_t^i=a^i\right)\cdot\beta^{t-1}\cdot\nonumber\\
	&\quad \sum_{\substack{s^j\in S_j,a^j\in A_j,\\j\in N,j\neq i}} w_i\left(\prod_{j\in N,j\neq i}\mathbb{P}^{\Psi^j}_{x^j}\left(X_t^j=s^j,A_t^j=a^j\right)\right) v_i\left(r_i(s^1,s^2,\cdots,s^n,a^1,a^2,\cdots,a^n)\right).\nonumber
\end{align}
Thus, $\max_{\pi^i\in\Pi_i}V_i^{\tilde{m}}(\boldsymbol{x},\pi^i,\boldsymbol{\Psi}^{-i})$ is an expected finite horizon Markov decision processes for player $i$ with nonstationary payoff function and transition probability. By \cite[Proposition 4.4.3]{puterman}, there exists a deterministic Markov policy that is optimal for $\max_{\pi^i\in\Pi_i}V_i^{\tilde{m}}(\boldsymbol{x},\pi^i,\boldsymbol{\Psi}^{-i})$. Thus, (\ref{finiteclose}) is equivalent to % to verify whether there exist a Markov strategy $\Psi\in\Pi_M$ satisfy (\ref{finiteclose}), we just need to consider
\begin{align}
	V_i^{\tilde{m}}(\boldsymbol{x},\Psi^i,\boldsymbol{\Psi}^{-i})\geq V_i^{\tilde{m}}(\boldsymbol{x},\pi^i,\boldsymbol{\Psi}^{-i})-\frac{\varepsilon}{3},~~\forall i\in N,\pi^i\in\Pi_{DM}^i.\nonumber
\end{align}
The advantage is that there are $|A_i|^{\tilde{m} \cdot |S_i|}$ elements (finite number) in $\Pi_{DM}^i$ for all $i\in N$ when considering the $\tilde{m}$-step finite truncation.

Below, we propose an iteration algorithm to find a Markov $\varepsilon$-equilibrium. 

\begin{algorithm}[H]\label{algo1}%[!htb]\label{algo1}
	\caption{An algorithm to solve a Markov $\varepsilon$-equilibrium of stochastic game (\ref{model})}
	{\bfseries Parameters:} Accuracy parameter $\varepsilon>0$, period length $\tilde{m}=\tilde{m}(\varepsilon)$ defined in Lemma \ref{finiteappro}, division parameter $\delta\in(0,1)$
%	the payoff rate $r(\Blue{t,}x,a,b)$, $\|r\|=\sup_{x\in S,\Blue{t\in [0,T],}a\in A\Blue{_t}(x),b\in B\Blue{_t}(x)}|r(\Blue{t,}x,a,b)|$; the transition rate $q(dy|\Blue{t,}x,a,b)$, $\|q\|=\sup_{x\in S} q^*(x)	=\sup_{x\in S,\Blue{t\in[0,T],}a\in A\Blue{_t}(x),b\in B\Blue{_t}(x)}|-q(\{x\}|\Blue{t,}x,a,b)|$; the terminal reward $g(x)$; the risk-sensitive parameter $\theta>0$;  finite horizon $T>0$; Player 1 has $|A|$ actions and Player 2	has $|B|$ actions, where $|A|$ is denoted the number of elements of a finite set $A$; a small error bound $\varepsilon>0$ determining the algorithm accuracy;	 the time step size after discretization is $d$%; contraction factor $\eta\gamma$ with $\varepsilon=\frac{\epsilon}{1-\eta\gamma}$; a measurable function $\omega:X\rightarrow [1,\infty)$
		\\
		{\bfseries Initialize:} Let $\boldsymbol{\Pi}^\delta_M:=\{\boldsymbol{\Psi}_1,\boldsymbol{\Psi}_2,\cdots,\boldsymbol{\Psi}_{k_{\delta}}\}$ be the set of all Markov strategies divided by the accuracy $\delta$ (i.e., for any $\boldsymbol{\Psi}\in\boldsymbol{\Pi}^\delta_M$, $\Psi^i_t(a^i|s^i)\in \{0,\delta,2\delta,\cdots,1\} \cup \{1,1-\delta,1-2\delta,\cdots,0\}$ holds for all $i\in N,t\in T,t\leq \tilde{m},s^i\in S_i,a^i\in A_i$); $\alpha=1$\\
%	$v_0(t,x)\in \mathbb{B}^{1}_{1,1}([0,T]\times S)$  arbitrarily, $n=0$\\
		{\bfseries Loop }\For{$\boldsymbol{\Psi}\in\boldsymbol{\Pi}^\delta_M$}
		{
			
			\For{$i=1:1:n$}
			{
				$\hat{V}_i \gets V_i^{\tilde{m}}(\boldsymbol{x},\Psi^i,\boldsymbol{\Psi}^{-i})$\\
				$\hat{V}_i^{max} \gets \max_{\pi^i\in \Pi_{DM}^i} V_i^{\tilde{m}}(\boldsymbol{x},\pi^i,\boldsymbol{\Psi}^{-i})$\\
				\uIf{$\hat{V}_i<\hat{V}_i^{max}-\frac{\varepsilon}{3}$}{
					\bfseries break
				}
				%\uElseIf{$i=n$}
				\ElseIf{$i=n$}
				{
					{\bfseries Output:}
					$\boldsymbol{\Psi}$, which is a Markov $\varepsilon$-equilibrium
				}
			%	\lElse{\bfseries continue}		
			}
			$\alpha \gets \alpha+1$\\
			\If{$\alpha=k_\delta$}{
				$\delta \gets \frac{\delta}{2}$\\
				Update the strategy set $\boldsymbol{\Pi}^\delta_M$ based on the accuracy $\delta$ (similar to the initialization)\\
				$\alpha \gets 1$
			}	
		}
	
%	\vspace*{-0.3cm}
\end{algorithm}

The convergence of Algorithm \ref{algo1} can be guaranteed, since the algorithm will stop within a finite number of cycles based on the above analysis.

\begin{remark}\label{delta_rem}
	\begin{enumerate}[(a)]
	\item If the weighting function $w_i~(i\in N)$ is Lipschitz continuous (i.e., there exists constant $C$ such that $|w_i(y)-w_i(z)|\leq C|y-z|$ for all $y,z\in [0,1]$)\footnote{For example, $w_i(y)=1-(1-y)^\alpha (\alpha > 1)$ for all $i\in N$ \cite[Example 8]{lin}.}, the division parameter $\delta=\delta(\varepsilon)$ satisfying (\ref{leqdelta}) and (\ref{aim}) can be given accurately as below. Without loss of generality, we assume that $C\geq 1$. Then, $\delta=\frac{(1-\beta)\varepsilon}{6CK|\boldsymbol{S}||\boldsymbol{A}|\hat{m}\sum_{i\in N}(|S_i||A_i|)^{\hat{m}}}$, which $K:= \max_{i\in N}\{K_i\},\hat{m} := \left\lceil \max_{i\in N}\left\{\frac{\ln{\frac{(1-\beta)\varepsilon}{6K_i |\boldsymbol{S}| |\boldsymbol{A}|}}}{\ln{\beta}}\right\} \right\rceil$.  The derivation is
	presented in Appendix \ref{appendix1}.
	\item When the division parameter $\delta$ is not analytically expressed, the strategy $\boldsymbol{\Psi}\in\boldsymbol{\Pi}_M$ can also be verified under the sampling without replacement method, instead of $\boldsymbol{\Psi}\in\boldsymbol{\Pi}^\delta_M$.
	\end{enumerate}
\end{remark}

\section{Applications in the Smart Grid}\label{applisection}

In this section, we present an application of nonzero-sum discounted decentralized stochastic games under PT.

\subsection{Discounted Stochastic Games among Prospect Prosumers}\label{prosusection}

On the supply side of the smart grid, a utility company sells energy to many prosumers who can produce energy from renewable resources and consume energy. On the demand side, prosumers are equipped with storage devices. Under an energy price rule given by the utility company, at each time slot, every prosumer decides the amount of energy they consume and purchase. Since the energy price rule is related to all prosumers' energy demand and the generated renewable energy is uncertain, prosumers interact with each other in a noncooperative stochastic game. In the followings, we investigate the interaction game among prosumers who exhibit subjective behavior, and focus on the time value of utility. Each prosumer strives to maximize their prospect utility.

As an application of model (\ref{model}), we introduce the nonzero-sum discounted decentralized stochastic games among $n$ prospect prosumers with renewable resources: 

\begin{enumerate}
	\item $N$ is prosumers set.
	\item State space $S_i=\{0,1,\cdots,\bar{S}_i\}$ is the space of energy storage for prosumer $i$ with maximum storage capacity $\bar{S}_i$. %$\boldsymbol{S}:=\prod_{i\in N}S_i$ is state space for the smart grid system and represents the space of energy storage. Assume that the energy storage in the initial time is $x=(x^1,x^2,\cdots,x^n)\in S$;
	\item $A_i=L_i\times D_i$ is action space for each prosumer $i$, where $L_i:=\{0,1,\cdots,\bar{L}_i\}$ (resp. $D_i:=\{0,1,\cdots,\bar{D}_i\}$) represents the space of energy consumption (resp. energy demand) for prosumer $i$ with maximum consumption $\bar{L}_i$ (maximum demand $\bar{D}_i$) at any time. %Denote $\boldsymbol{A}:=\prod_{i\in N}A_i$, i.e., action space of all prosumers;% $A_i(x)\subset A_i$ is the available action set for player $i\in I$ in state $x\in S$. 
	\item $r_i(\boldsymbol{s},\boldsymbol{a})$ is payoff function of prosumer $i$, where $\boldsymbol{s}=(s^1,s^2,\cdots,s^n)\in \boldsymbol{S}, \boldsymbol{a}=(l^1,d^1,l^2,d^2,\cdots,l^n,d^n)\in \boldsymbol{A}$. Define $\boldsymbol{d}:=(d^i:i\in N)$. Assume that
	\begin{align}
		r_i(\boldsymbol{s},\boldsymbol{a})=f_i(l^i)-d^i p_i(\boldsymbol{d})-g_i(s^i).\nonumber
	\end{align}
	Here, $f_i$ and $g_i$ are the satisfaction function for consuming $l^i$ energy and the storage cost function\footnote{The cost of energy storage system includes capacity cost (opportunity cost) and energy consumption cost \cite{liu_cost}. Here, $g_i$ primarily refers to the opportunity cost. If this cost is too low not to affect prosumers, let $g_i\equiv 0$.
	} for $s^i$ energy of prosumer $i$, respectively. $p_i$ is the energy price function, decided by the utility company according to all prosumers' demand.
	Actually, the payoff function is independent on the energy storage of other prosumers.%, we have $r_i(\boldsymbol{s},\boldsymbol{a})=:r_i(s^i,\boldsymbol{a})(\forall i\in N)$;
	\item 
	A nonstationary state transition equation is given by,
	\begin{align}
		X^i_{t+1}=\min\left\{\max\{0,X^i_t+G^i_t+D^i_t-L^i_t\},\bar{S}_i\right\},~~\forall i\in N,t\in T.\nonumber
	\end{align}
	For prosumer $i$ at time $t$, $X^i_t$ and $D^i_t$ are energy storage and energy demand respectively; $G^i_t$ represents the effective generated energy harvested from renewable resources minus the storage self-discharge\footnote{It can be explained as the effective generated energy accumulated during time $t_0\in(t-1,t]$, which prosumers do not know at time $t$. They can consume it after time $t$.}, which is a random variable on a probability space $(\Omega,\mathscr{F},\mathbb{P})$ with integer values. $L^i_t$ is energy consumption for prosumer $i$ during time $t_0\in(t,t+1]$.
	Then, the transition probability is determined by the nonstationary state transition equation as below
	\begin{align}
		q_t^i(y^i|s^i,a^i)&:=\mathbb{P}\left(X^i_{t+1}=y^i\big| X^i_{t}=s^i,A^i_{t}=a^i\right),\quad i\in N,\nonumber\\
		q_t(\boldsymbol{y}|\boldsymbol{s},\boldsymbol{a})
		&:=\mathbb{P}\left(X^i_{t+1}=y^i,i\in N \big| X^i_{t}=s^i,A^i_{t}=a^i,i\in N\right),\nonumber
	\end{align}
	for all $t\in T$, where $A^i_{t}:=(L^i_{t},D^i_{t}), a^i\in A_i$ and $s^i,y^i\in S_i$ ($\boldsymbol{a},\boldsymbol{s},\boldsymbol{y}$ are the vector of $a^i,s^i,y^i$ about $i\in N$). Thus, the distribution of $G^i_t$ determines the transition probability.% Under the independence assumption among prosumers, the transition probability is composed of the state transition probability with each prosumer, i.e., 
	%	\begin{align}\label{probabilityindependence}
	%		q(y|s,{a})=\prod_{i\in N}q_i(y^i|s^i,a^i),
	%	\end{align}
	%	for all $s=(s^i:i\in N),y=(y^i:i\in N)\in S,a=(a^i:i\in N)\in A$;
%	\item $\beta \in (0,1)$ is a discount factor.
\end{enumerate}

%Next, we give an informal description of the evolution of model (\ref{model}) as follows: In the initial state $x\in S$, $m$ players independently and simultaneously choose their actions {$a_i$} from $A_i(x)(i\in I)$ and obtain payoffs $r(i,x,\{a}){(\{a}=(a_1,a_2,\cdots,a_m))}$ after observing state $x$, and the system state jumps to $B\in \mathscr{B}(S)$ with the transition probability $q(B|x,\{a})$. Then this process repeats with payoffs discounted by $\beta$.

To verify the assumption of (\ref{probabilityindependence}), we propose some independence assumption among prosumers. %Before describing the purpose of the stochastic games, we hope to explore how the transition probability is determined by the distribution of $G^i_t$. First, we propose some independence assumption among prosumers.
\begin{assumption}\label{ass1}
	For any $t\in T$, the effective generated energy $G^1_t,G^2_t,\cdots,G^n_t$ are independent of each prosumer.
\end{assumption}

In the smart grid, prosumers are equipped with all kinds of renewable resource generators, such as solar, wind, tide, geothermal, etc. Whenever the generators between the prosumers are different, or if prosumers with the same generators are located at relatively far distances, the effective generated energy will be independent of each prosumer at every time. Thus, Assumption \ref{ass1} is reasonable.

\begin{remark}\label{remark1}
	Assumption \ref{ass1} is weaker than \cite[Assumption 1]{etesami} that studies the expected average criterion under PT. \cite{etesami} also assumes that for any $i\in N$:
	\begin{enumerate}[(a)]
		\item For any $t_1,t_2\in T$, $G^i_{t_1}$ and $G^i_{t_2}$ are identical distribution;
		\item $\min_{|k|\leq \bar{S}_i}\mathbb{P}(G^i_t=k)>0$ for all $t\in T$.
	\end{enumerate}
	The ergodicity under stationary strategy of each prosumer is important in analyzing the expected average criterion, which can be derived by (a) and (b). However, the expected discounted criterion does not need the ergodicity.
\end{remark}

%For any prosumer $i\in N$, define the state transition probability by
%\begin{align}
%	q_t^i(y^i|s^i,a^i):=&\mathbb{P}\left(X^i_{t+1}=y^i\big| X^i_{t}=s^i,A^i_{t}=a^i\right)\nonumber\\
%	=&\mathbb{P}\left(X^i_{t+1}=y^i\big| X^i_{t}=s^i,L^i_{t}=l^i,D^i_{t}=d^i\right)\nonumber\\
%	=&\left \{
%	\begin{array}{ll}
%		\mathbb{P}(G^i_t= l^i + y^i - s^i - d^i),      &  y^i = 1,2,\cdots,\bar{S}_i -1, \\
%		\sum_{\substack{k\in\mathbb{Z},\\s^i + k + d^i - l^i \leq 0}}\mathbb{P}(G^i_t=k),      & y^i = 0,\\
%		\sum_{\substack{k\in\mathbb{Z},\\s^i + k + d^i - l^i \geq \bar{S}_i}} \mathbb{P}(G^i_t=k),      &  y^i = \bar{S}_i,
%	\end{array}
%	\right.\nonumber
%\end{align}
%for all $t\in T$, where %$A^i_{t}=(L^i_{t},D^i_{t})$, 
%$a^i=(l^i,d^i)\in A_i$ and $s^i,y^i\in S_i$. Next, we describe how the transition probability is determined by the distribution of $G^i_t$, i.e., the relations of the state transition probability between the smart grid system and each prosumer.
\begin{proposition}\label{prop2}
	Under Assumption \ref{ass1}, the transition probability of the system is composed of the state transition probability with each prosumer as follows:
	\begin{align*}
		q_t(\boldsymbol{y}|\boldsymbol{s},\boldsymbol{a})=\prod_{i\in N}q_t^i(y^i|s^i,a^i),~~\forall t\in T,
	\end{align*}
for all $\boldsymbol{s}=(s^i:i\in N),\boldsymbol{y}=(y^i:i\in N)\in \boldsymbol{S},\boldsymbol{a}=(a^i:i\in N)\in \boldsymbol{A}$. Thus, (\ref{probabilityindependence}) holds.%Now, we denote $q_t(y|s,{a})(\forall t\in T)$ by $q(y|s,{a})$.
\end{proposition}
{\it Proof:} 
\begin{align}\label{prop1proof}
	&q_t(\boldsymbol{y}|\boldsymbol{s},\boldsymbol{a})\nonumber\\
	&=\mathbb{P}\left(X^i_{t+1}=y^i,i\in N \big| X^i_{t}=s^i,A^i_{t}=a^i,i\in N\right)\nonumber\\
	&=\mathbb{P}\left(X^1_{t+1}=y^1 \big| X^i_{t}=s^i,A^i_{t}=a^i,i\in N\right)\cdot\nonumber\\
	&~~~~~~
	\mathbb{P}\left(X^2_{t+1}=y^2,X^3_{t+1}=y^3,\cdots,X^n_{t+1}=y^n \big| X^i_{t}=s^i,A^i_{t}=a^i,i\in N,X^1_{t+1}=y^1\right)\nonumber\\
	&=\mathbb{P}\left(X^1_{t+1}=y^1 \big| X^1_{t}=s^1,A^1_{t}=a^1\right)\cdot\nonumber\\
	&~~~~~~
	\mathbb{P}\left(X^2_{t+1}=y^2,X^3_{t+1}=y^3,\cdots,X^n_{t+1}=y^n \big| X^i_{t}=s^i,A^i_{t}=a^i,i\in N,X^1_{t+1}=y^1\right)\nonumber\\
	&=\mathbb{P}\left(X^1_{t+1}=y^1 \big| X^1_{t}=s^1,A^1_{t}=a^1\right)\cdot\nonumber\\
	&~~~~~~\mathbb{P}\left(X^2_{t+1}=y^2 \big| X^i_{t}=s^i,A^i_{t}=a^i,i\in N,X^1_{t+1}=y^1\right)\cdot\nonumber\\
	&~~~~~~
	\mathbb{P}\left(X^3_{t+1}=y^3,X^4_{t+1}=y^4,\cdots,X^n_{t+1}=y^n \big| X^i_{t}=s^i,A^i_{t}=a^i,i\in N,X^1_{t+1}=y^1,X^2_{t+1}=y^2\right)\nonumber\\
	&=\mathbb{P}\left(X^1_{t+1}=y^1 \big| X^1_{t}=s^1,A^1_{t}=a^1\right)\cdot\nonumber\\
	&~~~~~~\mathbb{P}\left(X^2_{t+1}=y^2 \big| X^2_{t}=s^2,A^2_{t}=a^2\right)\cdot\nonumber\\
	&~~~~~~
	\mathbb{P}\left(X^3_{t+1}=y^3,X^4_{t+1}=y^4,\cdots,X^n_{t+1}=y^n \big| X^i_{t}=s^i,A^i_{t}=a^i,i\in N,X^1_{t+1}=y^1,X^2_{t+1}=y^2\right),
\end{align}
which the last equality is from Assumption \ref{ass1}. Then, similar to the derivation of (\ref{prop1proof}), we have
\begin{align}
	q_t(\boldsymbol{y}|\boldsymbol{s},\boldsymbol{a})=\prod_{i\in N}\mathbb{P}\left(X^i_{t+1}=y^i \big| X^i_{t}=s^i,A^i_{t}=a^i\right)=\prod_{i\in N}q_t^i(y^i|s^i,a^i),~~\forall t\in T.\nonumber
\end{align}
\te

From proposition \ref{prop2}, the assumption of (\ref{probabilityindependence}) is practical and meaningful in the smart grid. In Section \ref{applisection}, the followings are all based on Assumption \ref{ass1}. We do not require renewable energy generation to be identically distributed with respect to time, i.e., the state transition probability can be nonstationary, which makes it more realistic.

By Theorem \ref{theorem1}, there exists a Nash equilibrium for the nonstationary decentralized stochastic game among $n$ prosumers with the expected discounted payoff function under PT. Especially, for a given initial state $\boldsymbol{x}\in \boldsymbol{S}$, there exists a Markov $\boldsymbol{x}$-Nash equilibrium.

\subsection{Simulation Results}\label{simulsection}
For the discounted decentralized stochastic games among prospect prosumers in Section \ref{prosusection}, we give some numerical examples of Algorithm \ref{algo1} with the random sampling method (see Remark \ref{delta_rem} (b)) when accuracy parameter $\varepsilon=0.01$. In order to make the results more convincing, the partial settings of the stochastic game model are based on the reference \cite{etesami}.

We consider three prosumers in the smart grid system, i.e., $N=\{1,2,3\}$. Each prosumer has the same energy storage state space (energy consumption space) with three storage levels (two consumption levels, consumes energy or not), i.e., $S_i=\{0,1,2\}$ ($L_i=\{0,1\}$) for any $i=1,2,3$. 
Assume that $\tau_i\in S_i$ is a threshold constant which represents the extra energy that prosumer $i$ hopes to save for the next time stage, and let $\tau_1=1,\tau_2=0,\tau_3=1$.
For every time $t\in T$ (e.g., everyday), let the energy demand of prosumer $i\in N$ be $D^i_t=\max\{0,\tau_i+L^i_t-X^i_t\}$. Thus, the energy demand of each prosumer can be determined by the energy consumption, the energy storage and the threshold constant, then, the action variable reduces to the energy consumption.

Let the satisfaction function for consuming energy be $f_i(l^i)=\ln(1+l^i)~(\forall i\in N)$. Assume that $p_i(\boldsymbol{d})$ is a kind of fairness pricing function such as $p_i(\boldsymbol{d})=\frac{d^i}{\sum_{j=1}^{3}d^j}$. Now, we are not currently considering the cost of energy storage, i.e., $g_i\equiv 0$. Thus, the payoff function of prosumer $i$ is
\begin{align}
	r_i(\boldsymbol{s},\boldsymbol{l},\boldsymbol{d})=\ln(1+l^i)-d^i \frac{d^i}{\sum_{j=1}^{3}d^j},\nonumber
\end{align}
here, $\boldsymbol{s}=(s^1,s^2,s^3)\in \prod_{i=1}^{3}S_i, \boldsymbol{l}=(l^1,l^2,l^3)\in \prod_{i=1}^{3}L_i, \boldsymbol{d}=(d^1,d^2,d^3)$.

For simplicity, we consider that the effective generated energy $G^i_t$ is identical distribution with respect to $t$ for each prosumer $i$, except for Assumption \ref{ass1}. Denote $G^i_t~(t\in T)$ by $G^i$. Assume that $G^i$ confirms to Gaussian distribution, $G^1\sim N(0.5,2), G^2\sim N(0.5,1)$ and $G^3\sim N(1,1)$. That is to say, the state transition equation (probability) is stationary.

The weighting function and the valuation function for prosumer $i~(\forall i\in N)$ under PT are respectively $w_i(y)=e^{-(-\ln{y})^{0.8}}$ and
$v_i(y)=\left \{
\begin{array}{ll}
	y^{0.5},      &  y\geq 0 \\
	-(-y)^{0.3},      & y<0
\end{array}
\right.$.

Let the discount factor $\beta=0.001$, then, $\tilde{m}=\tilde{m}(\varepsilon)$ can be obtained from Lemma \ref{finiteappro}, we have $\tilde{m}=2$. The expected discounted payoff function under PT for prosumer $i\in N$ can be approximated by $2$-step finite truncation. Therefore, we mainly show the $2$-stage $\varepsilon$-equilibrium strategies in the followings.

%\begin{figure}[H]   
%	\centering         
%	\includegraphics[height=2cm]{Figure_1.pdf} 
%	\caption{}
%\end{figure}

\begin{figure}[t]%[htbp]	 
	\centering
	\captionsetup[subfigure]{labelformat=empty}
	\subfigure{
		\begin{minipage}[t]{0.33\linewidth}
			\centering
			\includegraphics[width=\linewidth,trim=3.8cm 1cm 2cm 0.8cm, clip]{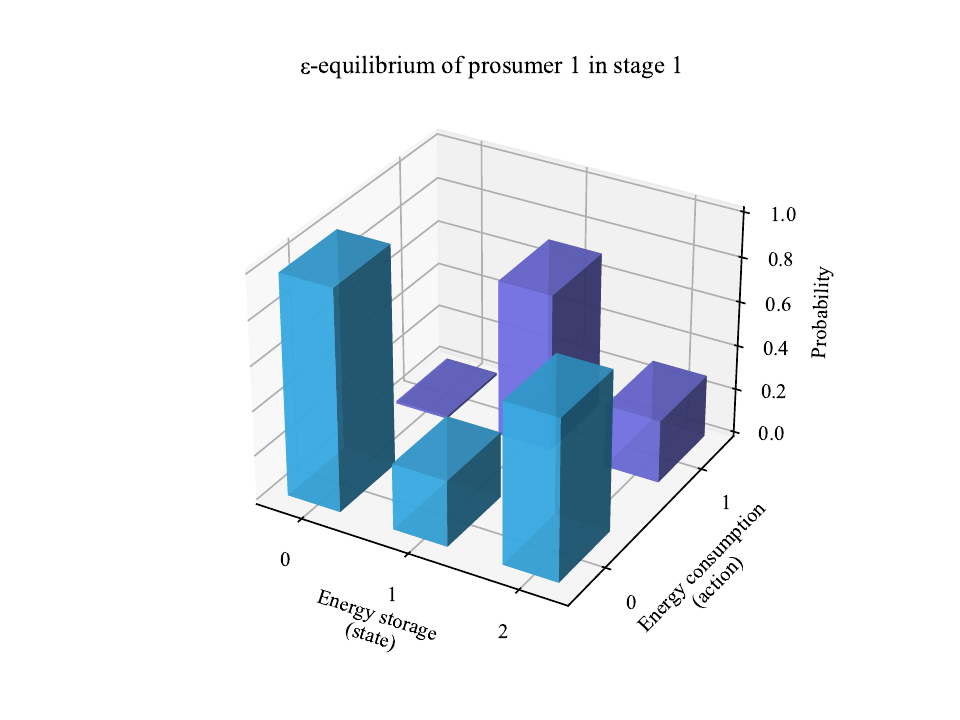}
			\caption*{}
		\end{minipage}%
		\begin{minipage}[t]{0.33\linewidth}
			\centering
			\includegraphics[width=\linewidth,trim=3.8cm 1cm 2cm 0.8cm, clip]{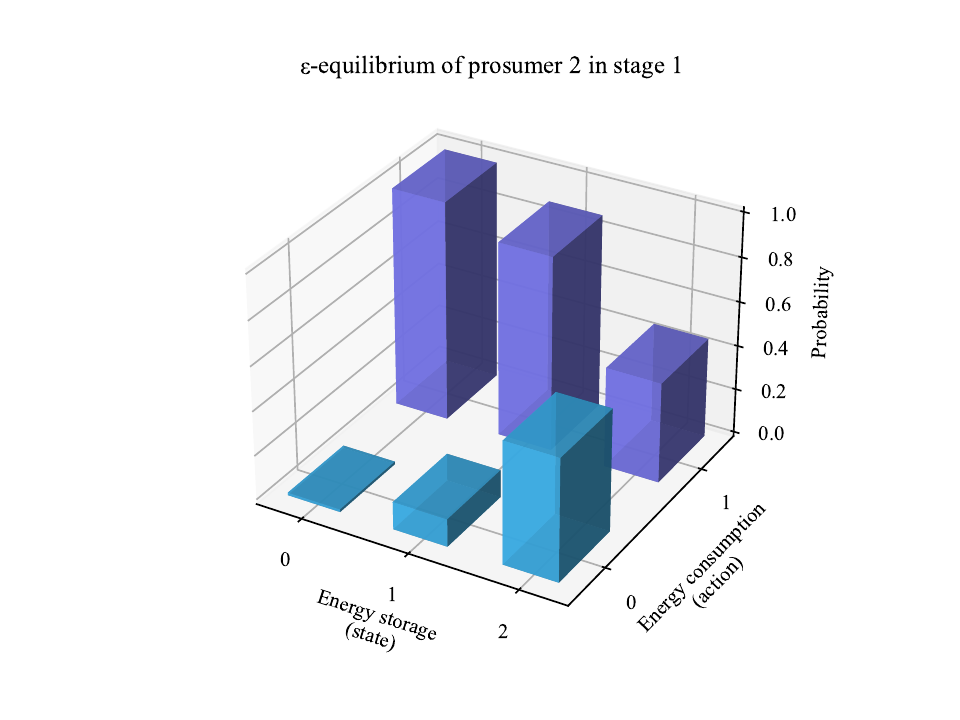}
			\caption*{}
		\end{minipage}%
		\begin{minipage}[t]{0.33\linewidth}
			\centering
			\includegraphics[width=\linewidth,trim=3.8cm 1cm 2cm 0.8cm, clip]{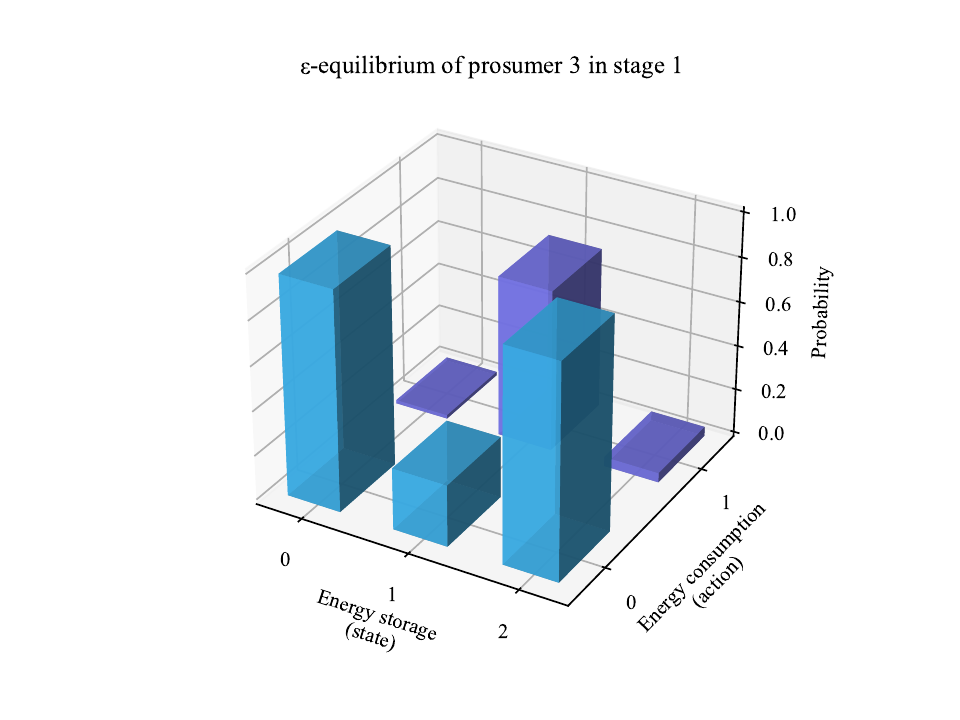}
			\caption*{}
	\end{minipage}}
\subfigure{
	\begin{minipage}[t]{0.33\linewidth} 
		\centering
		\includegraphics[width=\linewidth,trim=3.8cm 1cm 2cm 0.8cm, clip]{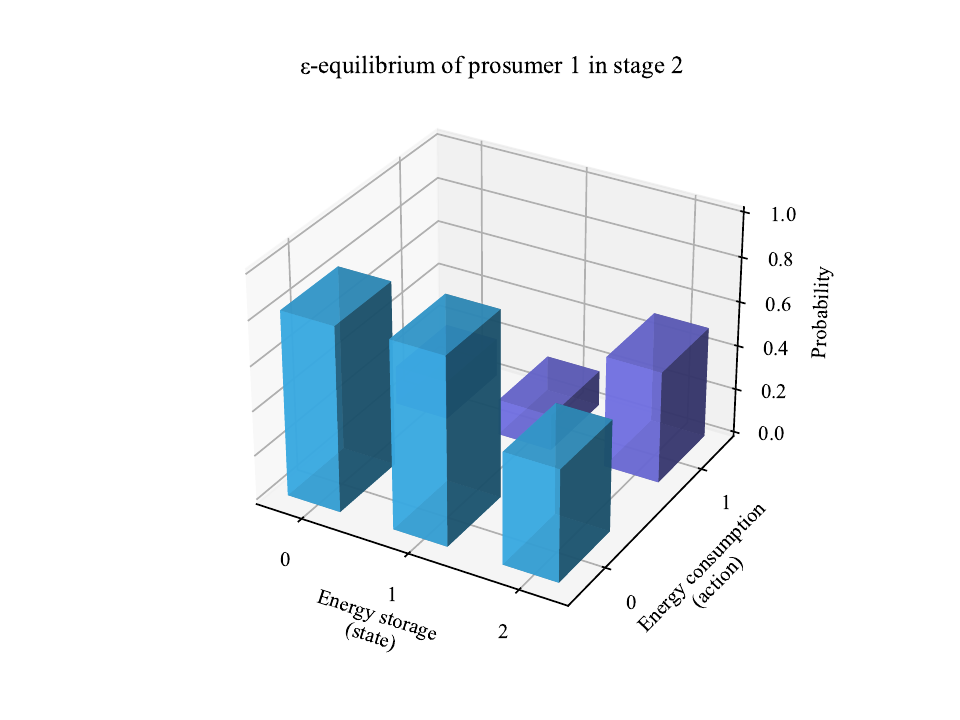} 
		\caption*{}
	\end{minipage}%
	\begin{minipage}[t]{0.33\linewidth}
		\centering
		\includegraphics[width=\linewidth,trim=3.8cm 1cm 2cm 0.8cm, clip]{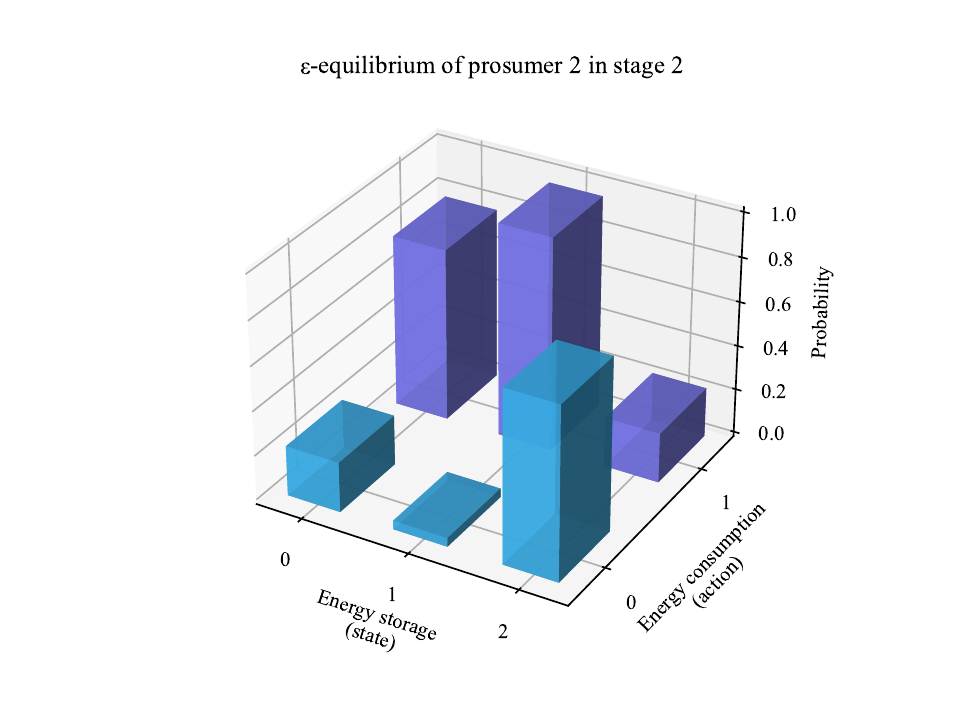}
		\caption*{}
	\end{minipage}%
	\begin{minipage}[t]{0.33\linewidth}
		\centering
		\includegraphics[width=\linewidth,trim=3.8cm 1cm 2cm 0.8cm, clip]{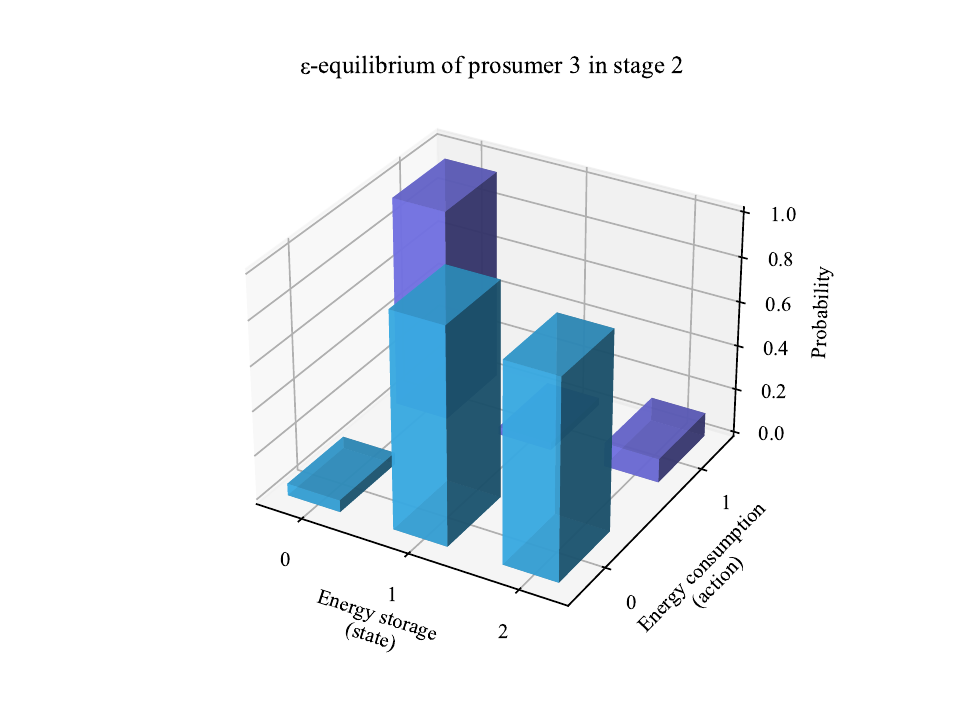}
		\caption*{}
\end{minipage}}
	\caption{Markov $\varepsilon$-equilibrium for initial energy storage state $\boldsymbol{x}=(0,0,0)$}
	\label{ENE000}
\end{figure}

For initial energy storage state $\boldsymbol{x}=(0,0,0)$, Figure \ref{ENE000} illustrates a Markov $\varepsilon$-equilibrium of all prosumers\footnote{Note that there may be multiple $\varepsilon$-equilibria.}, in which the height of the bar in $(s,l)$th coordinate represents the probability for consuming $l$ energy when this prosumer has $s$ energy storage. In Figure \ref{ENE000}, the followings can be seen.
\begin{enumerate}
	\item The equilibrium of the expected discounted criterion under PT is indeed nonstationary, though the stochastic game model is stationary. It is different from the expected average criterion under PT \cite{etesami} which finds a stationary $\varepsilon$-equilibrium.
	\item The equilibrium behaviors of prosumer 1 and 3 are similar in stage 1. The second stage reveals distinction, since the different normal distribution of the effective generated energy $G^1,G^3$. The mean of the generated energy for prosumer 3 is higher, so that prosumer 3 tends to consume without energy storage.
	\item The equilibrium behavior of prosumer 2 exhibits significant differences from prosumer 1 and 3. Prosumer 2 is more likely to consume energy at a lower energy storage level. This is because the threshold constant $\tau_i$ is different, and prosumer 2 does not need to save extra energy for the next time stage. 
	\item In a higher energy storage level, all prosumers are more inclined to not consume energy, or to balance whether energy is used. Consuming energy brings greater satisfaction. But in the while, it leads to the purchase of more energy from the utility company and increases the cost. Furthermore, it is consistent with prospect theory that prosumers are more conservative for higher energy storage due to the uncertainties.
\end{enumerate}

%\Red{initial energy storage state $\boldsymbol{x}=(1,1,1)$ and $\boldsymbol{x}=(2,2,2)$}

\section{Conclusions}\label{conclusection}

In this paper, we are concerned with the nonzero-sum decentralized stochastic games and applications in noncooperative games among prosumers on the demand side of the smart grid. According to prospect theory, a novel optimization objective function is the expected discounted payoff function, which characterizes the time value of utility and the subjective behavior of players. The state transition probability and payoff function can be nonstationary, which is more general. Significantly, several previous tools, such as the optimality equation, the occupation measure, and linear programming, are not feasible. We demonstrate the existence of Nash equilibrium by a new technique for constructing the marginal distribution on the state-action pairs at any time. Moreover, we develop an algorithm for $\varepsilon$-equilibrium and discuss the simulation results.

When confronted with a substantial scale of stochastic game models, the curse of dimensionality of Algorithm \ref{algo1} arises. As a future work of research, several learning algorithms, such as reinforcement learning, can be considered. In the application of the smart grid, it is possible to incorporate the utility company as a participant, based on the stochastic games among prosumers.

%\begin{appendices}
%\section{Derivation of the division parameter $\delta$}\label{appendix1}
%\end{appendices}

\appendix

\section{Derivation of the division parameter $\delta$ in Remark \ref{delta_rem} (a)}\label{appendix1}
Similar to Lemma \ref{finiteappro}, there exists a constant $\hat{m}:= \left\lceil \max_{i\in N}\left\{\frac{\ln{\frac{(1-\beta)\varepsilon}{6K_i |\boldsymbol{S}| |\boldsymbol{A}|}}}{\ln{\beta}}\right\} \right\rceil$ such that
\begin{align}
	\bigg|V_i(\boldsymbol{x},\boldsymbol{\phi})-V_i^{\hat{m}}(\boldsymbol{x},\boldsymbol{\phi})\bigg|\leq \frac{\varepsilon}{6},~~\forall i\in N, \boldsymbol{\phi}\in\boldsymbol{\Pi}.\nonumber
\end{align}

For two Markov strategies $\boldsymbol{\pi},\boldsymbol{\psi}\in\boldsymbol{\Pi}_M$, note that 
\begin{align}
	&|V_i(\boldsymbol{x},\boldsymbol{\pi})-V_i(\boldsymbol{x},\boldsymbol{\psi})|\nonumber\\
	&\leq |V_i(\boldsymbol{x},\boldsymbol{\pi})-V_i^{\hat{m}}(\boldsymbol{x},\boldsymbol{\pi})| + |V_i^{\hat{m}}(\boldsymbol{x},\boldsymbol{\pi})-V_i^{\hat{m}}(\boldsymbol{x},\boldsymbol{\psi})| + |V_i^{\hat{m}}(\boldsymbol{x},\boldsymbol{\psi})-V_i(\boldsymbol{x},\boldsymbol{\psi})|\nonumber\\
	&\leq \frac{\varepsilon}{3} + |V_i^{\hat{m}}(\boldsymbol{x},\boldsymbol{\pi})-V_i^{\hat{m}}(\boldsymbol{x},\boldsymbol{\psi})|,\nonumber
\end{align}
and we hope to obtain $|V_i^{\hat{m}}(\boldsymbol{x},\boldsymbol{\pi})-V_i^{\hat{m}}(\boldsymbol{x},\boldsymbol{\psi})|\leq \frac{\varepsilon}{6}$ for all $i\in N$, then $|V_i(\boldsymbol{x},\boldsymbol{\pi})-V_i(\boldsymbol{x},\boldsymbol{\psi})|\leq \frac{\varepsilon}{2}$. For any $i\in N,s^i\in S_i,a^i\in A_i,t\in T$, consider 
\begin{align}
	&\left| \mathbb{P}^{\pi^i}_{x^i}\left(X_t^i=s^i,A_t^i=a^i\right) - \mathbb{P}^{\psi^{i}}_{x^i}\left(X_t^i=s^i,A_t^i=a^i\right) \right| \nonumber\\
	&=\left|\sum_{\substack{s^i_2,s^i_3,\cdots,s^i_{t-1}\in S_i,\\a^i_1,a^i_2,\cdots,a^i_{t-1}\in A_i}}\pi^i_1(a^i_1|x^i)q^i_1(s^i_2|x^i,a^i_1)\pi^i_2(a^i_2|s^i_2)\cdots q^i_{t-1}(s^i|s^i_{t-1},a^i_{t-1})\pi^{i}_t(a^i|s^i)\right.\nonumber\\
	&\quad \left.-\sum_{\substack{s^i_2,s^i_3,\cdots,s^i_{t-1}\in S_i,\\a^i_1,a^i_2,\cdots,a^i_{t-1}\in A_i}}\psi^{i}_1(a^i_1|x^i)q^i_1(s^i_2|x^i,a^i_1)\psi^{i}_2(a^i_2|s^i_2)\cdots q^i_{t-1}(s^i|s^i_{t-1},a^i_{t-1})\psi^{i}_t(a^i|s^i)\right|\nonumber\\
	&\leq \sum_{\substack{s^i_2,s^i_3,\cdots,s^i_{t-1}\in S_i,\\a^i_1,a^i_2,\cdots,a^i_{t-1}\in A_i}}q^i_1(s^i_2|x^i,a^i_1)q^i_2(s^i_3|s^i_2,a^i_2)\cdots q^i_{t-1}(s^i|s^i_{t-1},a^i_{t-1})\left|\pi^{i}_1(a^i_1|x^i)\pi^{i}_2(a^i_2|s^i_2)\cdots \right.\nonumber\\
	&~~~~~~~~~~~~~~~~~~~~~~~~~~~~~~\left.\pi^{i}_t(a^i|s^i) - \psi^{i}_1(a^i_1|x^i)\psi^{i}_2(a^i_2|s^i_2)\cdots \psi^{i}_t(a^i|s^i) \right|\nonumber\\
	&\leq \sum_{\substack{s^i_2,s^i_3,\cdots,s^i_{t-1}\in S_i,\\a^i_1,a^i_2,\cdots,a^i_{t-1}\in A_i}}\left|\pi^{i}_1(a^i_1|x^i)\pi^{i}_2(a^i_2|s^i_2)\cdots \pi^{i}_t(a^i|s^i) - \psi^{i}_1(a^i_1|x^i)\psi^{i}_2(a^i_2|s^i_2)\cdots \psi^{i}_t(a^i|s^i) \right|\nonumber\\
	&\leq t\delta |S_i|^{t-2} |A_i|^{t-1}.\nonumber
\end{align}
And, 
\begin{align}
	\left| \prod_{\substack{j\in N\\j\neq i}}\mathbb{P}^{\pi^j}_{x^j}\left(X_t^j=s^j,A_t^j=a^j\right) - \prod_{\substack{j\in N\\j\neq i}}\mathbb{P}^{\psi^j}_{x^j}\left(X_t^j=s^j,A_t^j=a^j\right) \right|\leq t\delta \sum_{j\in N,j\neq i}|S_j|^{t-2} |A_j|^{t-1}.\nonumber
\end{align}
Then, by the Lipschitz continuity of $w_i~(i\in N)$, 
\begin{small}
\begin{align}\label{derivation}
	&|V_i^{\hat{m}}(\boldsymbol{x},\boldsymbol{\pi})-V_i^{\hat{m}}(\boldsymbol{x},\boldsymbol{\psi})| \nonumber\\
	&= \left|\sum_{t=1}^{\hat{m}}\beta^{t-1}\sum_{\boldsymbol{s}\in \boldsymbol{S},\boldsymbol{a}\in \boldsymbol{A}}v_i\left(r_i(\boldsymbol{s},\boldsymbol{a})\right)\left[\mathbb{P}_{i,\rm{PT}}^{\boldsymbol{x},\boldsymbol{\pi}}\left(\tilde{\boldsymbol{X}}_t=\boldsymbol{s},\tilde{\boldsymbol{A}}_t=\boldsymbol{a}\right) - \mathbb{P}_{i,\rm{PT}}^{\boldsymbol{x},\boldsymbol{\psi}}\left(\tilde{\boldsymbol{X}}_t=\boldsymbol{s},\tilde{\boldsymbol{A}}_t=\boldsymbol{a}\right)\right]\right|\nonumber\\
	&=\left|\sum_{t=1}^{\hat{m}}\beta^{t-1}\sum_{\boldsymbol{s}\in \boldsymbol{S},\boldsymbol{a}\in \boldsymbol{A}}v_i\left(r_i(\boldsymbol{s},\boldsymbol{a})\right)\left[	\mathbb{P}^{\pi^i}_{x^i}\left(X_t^i=s^i,A_t^i=a^i\right)w_i\left(\prod_{j\in N,j\neq i}\mathbb{P}^{\pi^j}_{x^j}\left(X_t^j=s^j,A_t^j=a^j\right)\right) \right.\right. \nonumber\\
	&\quad - \left.\left.	\mathbb{P}^{\psi^i}_{x^i}\left(X_t^i=s^i,A_t^i=a^i\right)w_i\left(\prod_{j\in N,j\neq i}\mathbb{P}^{\psi^j}_{x^j}\left(X_t^j=s^j,A_t^j=a^j\right)\right)\right]\right|\nonumber\\
	&\leq \sum_{t=1}^{\hat{m}}\beta^{t-1}\sum_{\boldsymbol{s}\in \boldsymbol{S},\boldsymbol{a}\in \boldsymbol{A}}v_i\left(r_i(\boldsymbol{s},\boldsymbol{a})\right) \cdot \nonumber\\
	&\quad \left[
	\left| w_i\left(\prod_{j\in N,j\neq i}\mathbb{P}^{\pi^j}_{x^j}\left(X_t^j=s^j,A_t^j=a^j\right)\right) - w_i\left(\prod_{j\in N,j\neq i}\mathbb{P}^{\psi^j}_{x^j}\left(X_t^j=s^j,A_t^j=a^j\right)\right) \right|  \right.\nonumber\\
	&\quad+ \left| \mathbb{P}^{\pi^i}_{x^i}\left(X_t^i=s^i,A_t^i=a^i\right) - \mathbb{P}^{\psi^i}_{x^i}\left(X_t^i=s^i,A_t^i=a^i\right) \right|
	\Bigg]\nonumber\\
	&\leq \sum_{t=1}^{\hat{m}}\beta^{t-1}K_i|\boldsymbol{S}||\boldsymbol{A}|\cdot \left( C\hat{m}\delta \sum_{j\in N,j\neq i}|S_j|^{\hat{m}-2} |A_j|^{\hat{m}-1} + \hat{m}\delta |S_i|^{\hat{m}-2} |A_i|^{\hat{m}-1} \right)\nonumber\\
	&\leq K|\boldsymbol{S}||\boldsymbol{A}|\cdot \left( C\hat{m}\delta \sum_{i\in N}(|S_i| |A_i|)^{\hat{m}} \right)\cdot \frac{1}{1-\beta}.
\end{align}
\end{small}
Let the right hand side of (\ref{derivation}) be $\frac{\varepsilon}{6}$, thus, $\delta=\frac{(1-\beta)\varepsilon}{6CK|\boldsymbol{S}||\boldsymbol{A}|\hat{m}\sum_{i\in N}(|S_i||A_i|)^{\hat{m}}}$.

%\vspace*{0.1in}

%\newpage
%\vspace{1in}

\end{document}